\documentclass{amsart}
\usepackage{graphicx}
\usepackage{psfrag}
\usepackage{amssymb}
\numberwithin{equation}{section}
\let\cal\mathcal

\def\Hscr{{\cal H}}

\let\blb\mathbb
\def\CC{{\blb C}}

\def\r{\rightarrow}

\newtheorem{lemma}{Lemma}[section]
\newtheorem{proposition}[lemma]{Proposition}

\theoremstyle{definition}

{

}

\theoremstyle{remark}

\newdimen\uboxsep \uboxsep=1ex
\def\uboxn#1{\vtop to 0pt{\hrule height 0pt depth 0pt\vskip\uboxsep
\hbox to 0pt{\hss #1\hss}\vss}}

\def\uboxs#1{\vbox to 0pt{\vss\hbox to 0pt{\hss #1\hss}
\vskip\uboxsep\hrule height 0pt depth 0pt}}

\def\Re{\operatorname{Re}}
\def\Im{\operatorname{Im}}
\keywords{Kontsevich weight}
\author{Michel Van den Bergh}
\address{Departement WNI\\Universiteit Hasselt\\ Universitaire Campus\\ Building D\\ 3590 
Diepenbeek\\ Belgium} 
\thanks{The author is a director of research at the FWO} 
\email{michel.vandenbergh@uhasselt.be} 
\title{The Kontsevich weight of a wheel with spokes pointing outward}
\dedicatory{Dedicated to  Fred Van Oystaeyen on the occasion of
his 60th birthday.}

\begin{document}
\begin{abstract}
  This is a companion note to ``Hochschild cohomology and Atiyah
  classes'' by Damien Calaque and the author. Using elementary methods
  we compute the Kontsevich weight of a wheel with spokes pointing
  outward. The result is in terms of modified Bernoulli numbers.
The same result had been obtained earlier by Torossian (unpublished) and
also recently by Thomas Willwacher using  more advanced methods.
\end{abstract}
\maketitle
\setcounter{tocdepth}{1}
\tableofcontents
\section{Introduction and statement of the main result}
\label{ref-1-0}
This is a companion note to \cite{vdbcalaque}. We refer to \cite[\S 8]{vdbcalaque} and \cite{Ko3} for unexplained notations and conventions. 

Below we will compute the Kontsevich weight $w_n$
of the ``opposite'' wheel
\[
\psfrag{m}[][]{}
\psfrag{w}[][]{$w$}
\psfrag{zn}[][]{$z_n$}
\psfrag{z1}[][]{$z_1$}
\psfrag{z2}[][]{$z_2$}
\includegraphics[width=4cm]{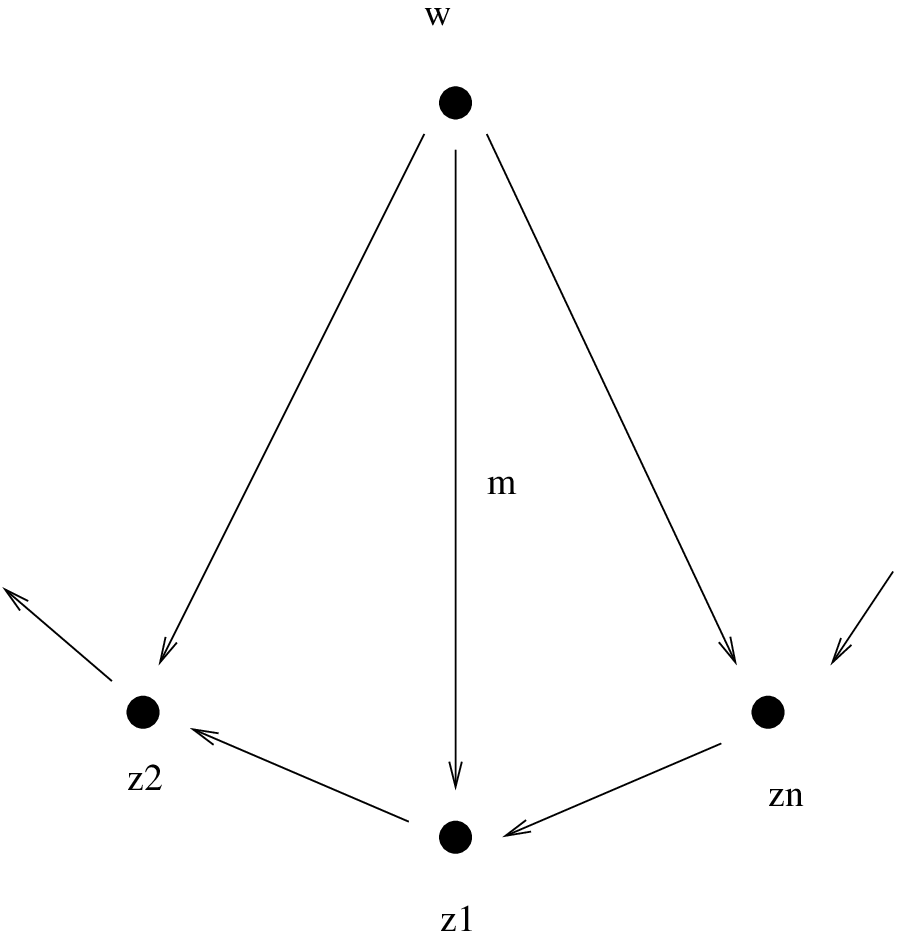}
\]
with edge ordering $(z_1,z_2)<\cdots<(z_n,z_1)<(w,z_1)<\cdots<(w,z_n)$. 

Thus\footnote{This is the traditional Kontsevich weight. In
  \cite[\S8]{vdbcalaque} we use $W_n=
  (-1)^{2n(2n-1)/2}w_n=(-1)^nw_n$. Since $w_n=0$ for $n$ odd we actually
have $W_n=w_n$.}
\[
w_n=\frac{1}{(2\pi)^{2n}}
\int_{C_{n+1,0}} d\varphi(z_1,z_2)\cdots d\varphi(z_1,z_n) d\varphi(w,z_1)\cdots d\varphi(w,z_n)
\]
where $\varphi(u,v)$ is computed as in the following diagram. 
\begin{center}
\psfrag{h}[l][l]{$\varphi(u,v)$}
\psfrag{p}[][]{$u$}
\psfrag{q}[][]{$v$}
\psfrag{e}[][]{$$}
\includegraphics[width=7cm]{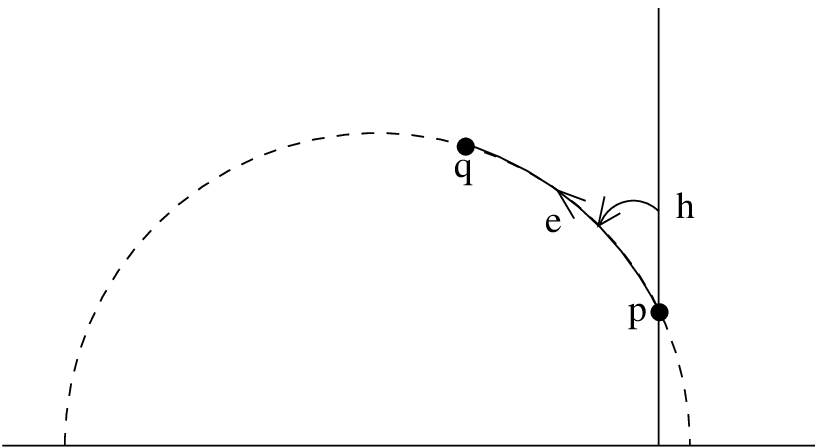}
\end{center}
and where $C_{n+1,0}$ is the configuration space of $n+1$ points in
the complex upper half plane $\Hscr$.  If we normalize things by
putting $w=i$ then the integration domain is equal to the complement of
the big diagonal
in $(\Hscr-\{i\})^n$.

Recall (e.g.\
\cite{Wolfram}) that the modified Bernoulli numbers $\ss_n$
are defined by
\begin{align*}
\sum_{n} \ss_n x^n&=\frac{1}{2}\log \frac{e^{x/2}-e^{-x/2}}{x}\\
&=\frac{x^2}{48}
-\frac{x^4}{5760}+\frac{x^6}{ 362880}-\frac{x^8}{19353600}+
\frac{x^{10}}{958003200} - \frac{691 x^{12}}{31384184832000}
+
\cdots
\end{align*}
The main result of this note is a relatively elementary proof of the following
(see \S\ref{ref-9.6-62})
\begin{proposition} \cite[Theorem 18]{Willwacher}
\begin{equation}
\label{ref-1.1-1}
w_n=-(-1)^{n(n-1)/2}n\ss_n
\end{equation}
\end{proposition}
We needed \eqref{ref-1.1-1} for the proof of \cite[Thm
1.3]{vdbcalaque}. At that time we didn't have a reference so we did
the computation ourselves. After this note was posted on the arXiv
Torossian sent us an argument which shows that \eqref{ref-1.1-1}
follows from the work of Cattaneo and Felder \cite{CF2,CF3} on the
quantization of coisotropic manifolds combined with Shoikhet's
vanishing results for wheels with spokes pointing inwards
\cite{Shoikhet} and some facts related to the
Campbell-Baker-Hausdorff formula. In this way one obtains a shorter and more
conceptual proof of \eqref{ref-1.1-1}.
Torossian also pointed out Willwacher's
recent work \cite{Willwacher} which uses the same ingredients. 

While more conceptual it seems this approach is substantially less
elementary than ours so we think that this note is still useful.

Let us finally also mention that a result similar to
\eqref{ref-1.1-1} was announced by Shoikhet
in \cite[\S2.3.1]{Shoikhet}.  It can presumably be obtained from the
methods in \cite{Shoikhet1}.

\medskip

The structure of this note is as follows: in \S\ref{ref-3-3} we represent
our computations graphically. This makes it easy to see
which boundary components in Stokes theorem yield non-trivial
contributions. We ignore signs as these are not easy to handle graphically.

In \S\ref{ref-4-14}-\S\ref{ref-8-47} we compute the contributions of the
relevant boundary components with precise signs. Thanks to some
``lucky'' cancellations we obtain a system of recursions for computing
$w_n$.

Finally in \S\ref{ref-9-53} we solve the recursions.

The reader will notice that  some formulas
and figures appear more than once in this paper. 
This has been done to make the three logical parts
\S\ref{ref-3-3},\S\ref{ref-4-14}-\S\ref{ref-8-47},\S\ref{ref-9-53} more or less
independently readable. This should make our computations easier to follow. 

\medskip

As indicated above we use almost nothing but Stokes theorem in our
computations. For a general introduction to the use of Stokes theorem in the
computation of Kontsevich weights see Torossian's article in \cite{CKTB}

\section{Acknowledgment}
\label{ref-2-2}
The symbolic computations were done using the SAGE \cite{SAGE}
interface to the computer algebra package Maxima \cite{Maxima}. Many
identities and integrals were tested numerically
using the Monte Carlo integration library which is part of GSL
\cite{GSL}.

We thank Charles Torossian for pointing how \eqref{ref-1.1-1} follows
from the work of Cattaneo and Felder.

Finally I would like to specially thank the anonymous referee for his very careful
reading of the manuscript. This has allowed me to improve the 
exposition and to correct numerous ``misprints''. 
\section{Graphical representation of  the computation}
\label{ref-3-3}
In this section we give a graphical overview of our computations. We will
ignore signs.  Hence we don't care about orientations or the ordering 
of the factors in a product of forms. 

\subsection{Graphical language}
\label{ref-3.1-4}
In the Kontsevich integral we associate to a graph $\Gamma$ a form $\omega_\Gamma$
on a suitable configuration space $C^+_{p,q}$. The form $\omega_\Gamma$ is the product of $d\varphi(u,v)$, where $(u,v)$ runs over the edges of $\Gamma$, divided by $(2\pi)^e$ where $e$ is the number of edges.

Besides the usual Kontsevich graphs we introduce graphs with some new
kind of arrows. 

An arrow of the form 
\[
\psfrag{n}[b][b]{$n$}
\psfrag{e}[t][t]{$$}
\psfrag{p}[t][t]{$u$}
\psfrag{q}[t][t]{$v$}
\includegraphics[width=2cm]{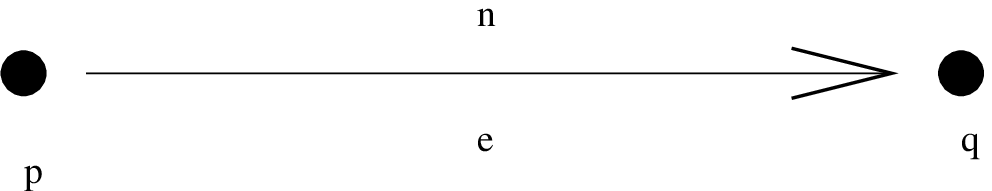}
\]
corresponds to a factor $d(\varphi(u,v)^n)/(2\pi)^n$ in $\omega_\Gamma$. 
For such an arrow
we assume that we have made a branch cut such that $\varphi(u,v)\in ]0,2\pi[$. 
I.e.\ if $\Re u=\Re v$ then $\Im v<\Im u$. Note that we have to specify the branch cut
as the value of the integral depends on it.

A bold arrow of the form
\[
\psfrag{n}[b][b]{$n$}
\psfrag{e}[t][t]{$$}
\psfrag{p}[t][t]{$u$}
\psfrag{q}[t][t]{$v$}
\includegraphics[width=2cm]{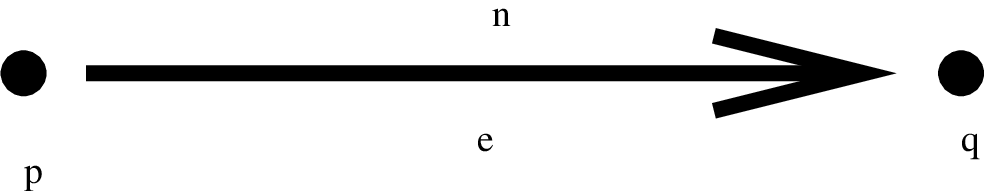}
\]
corresponds to a factor $\varphi(u,v)^n/(2\pi)^n$  in $\omega_\Gamma$. We make the same
branch cut as above.  If $n$ is not indicated then we assume $n=1$.

We introduce a third kind of special arrow. A dashed arrow of the form
\[
\psfrag{e}[t][t]{$$}
\psfrag{p}[t][t]{$u$}
\psfrag{q}[t][t]{$v$}
\includegraphics[width=2cm]{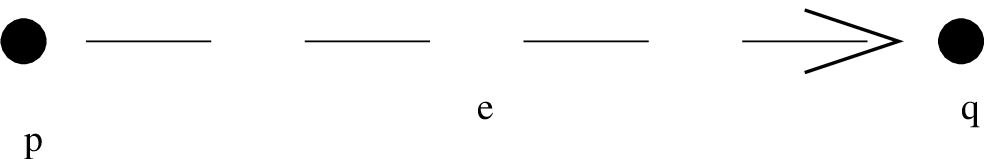}
\]
indicates a restriction of the integration domain to $\Re v=\Re u+\epsilon$,
$\Im v>\Im u$ with $\epsilon$ a positive infinitesimal. 

Here is the basic computation rule for computing $\int \omega_\Gamma$ for one of the
generalized graphs (Rule 1).
{\psfrag{n}[b][b]{$n$}
\psfrag{e}[t][t]{$$}
\psfrag{p}[t][t]{$u$}
\psfrag{q}[t][t]{$v$}
\begin{align*}
\int \includegraphics[width=2cm]{arrow_with_multiplicity.eps}&=
\pm\int_\partial \includegraphics[width=2cm]{bold_arrow.eps}
\end{align*}
\begin{align*}
\int_\partial \includegraphics[width=2cm]{bold_arrow.eps}
&=\left(\int_{\partial_{\text{classical}}} \includegraphics[width=2cm]{bold_arrow.eps}\right)\pm \left(\int \includegraphics[width=2cm]{dotted_arrow.eps}\right)
\end{align*}
}
Here $\partial$ denotes the boundary of the integration domain (including the
branch cut). $\partial_{\text{classical}}$ denotes the boundary without
the branch cut. $\partial_{\text{classical}}$ can be determined combinatorially
by contracting vertices in $\Gamma$ and moving vertices to the real line. 

The difference between $\partial$ and $\partial_{\text{classical}}$ is given by the
two sides of the branch cut. The difference in value of $\varphi(u,v)$ on
opposite sides of the branch cut is $2\pi$. This gives the contribution with the
dashed arrow. 

Here is another obvious computation rule (Rule 2)
\[
\psfrag{n}[b][b]{$n$}
\psfrag{n1}[b][b]{$n+1$}
\psfrag{e}[t][t]{$e$}
\psfrag{p}[t][t]{$p$}
\psfrag{q}[t][t]{$q$}
\raisebox{-2.1mm}{\includegraphics[width=2cm]{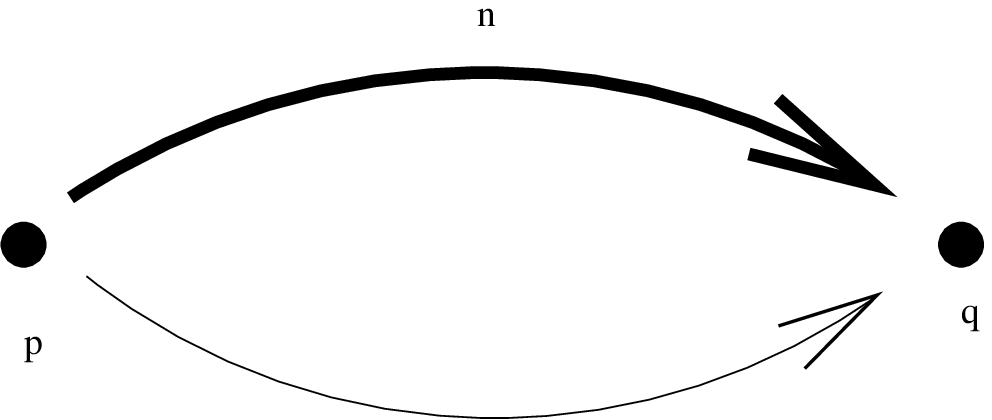}}\quad=\quad
\frac{1}{n+1}\times %\left(
\raisebox{-1.5mm}{\includegraphics[width=2cm]{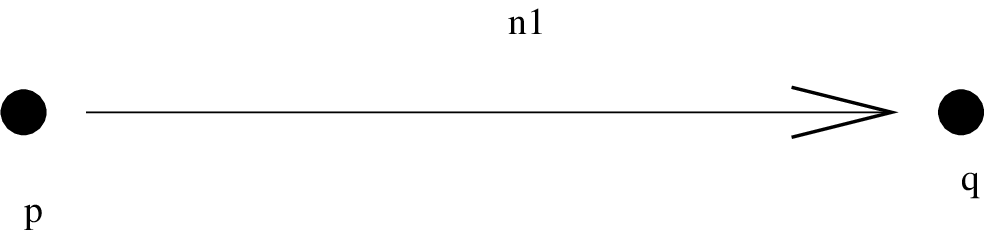}}
%\right)
\]
Below we will also consider some integrals over the configuration
spaces $C_p$ consisting of $p$ points in the complex plane. These are it terms of $\theta(u,v)$ where $\theta(u,v)$ 
is as in the following diagram 
\begin{center}
\psfrag{h}[l][l]{$\theta(u,v)$}
\psfrag{p}[][]{$u$}
\psfrag{q}[][]{$v$}
\psfrag{e}[][]{$$}
\includegraphics[width=2cm]{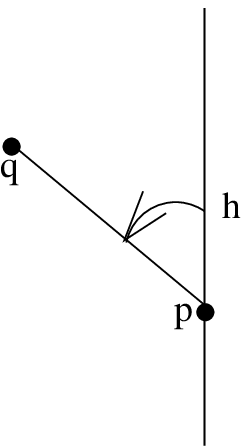}
\end{center}
To indicate that we are integrating over some $C_p$ we will put the graph
in a box. Otherwise we use the same conventions as above (with
$\theta(u,v)$ replacing $\varphi(u,v)$). Since $d\theta(u,v)=d\theta(v,u)$
some arrow indicators are superfluous. We will omit those. 

We will repeatedly use the following vanishing result
\begin{equation}
\label{lem66}
\int \raisebox{-0.4cm}{\includegraphics[width=1cm]{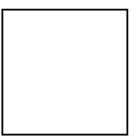}}=0
\end{equation}
provided the graph in the box contains only ordinary edges and
the number of edges is $\ge 3$. An integral as in \eqref{lem66} often arises
as a factor of an integral over a boundary component of $C_{n+1,0}$. 
\subsection{Applying Stokes theorem to $w_n$}
\label{ref-3.2-5}
We will now apply this graphical language to computations with wheels.
We assume first $n\ge 2$ 

The first step is (using Rule 1)
\begin{equation}
\label{ref-3.1-6}
\psfrag{w}[][]{$w$}
\psfrag{zn}[][]{$z_n$}
\psfrag{z1}[][]{$z_1$}
\psfrag{z2}[][]{$z_2$}
\psfrag{m}[l][l]{$$}
\int_{C_{n+1,0}} \raisebox{-2.5cm}{\includegraphics[width=4cm]{wheel_n.eps}}=\pm
\int_{\partial_{\text{classical}}}\raisebox{-2.5cm}{ \includegraphics[width=4cm]{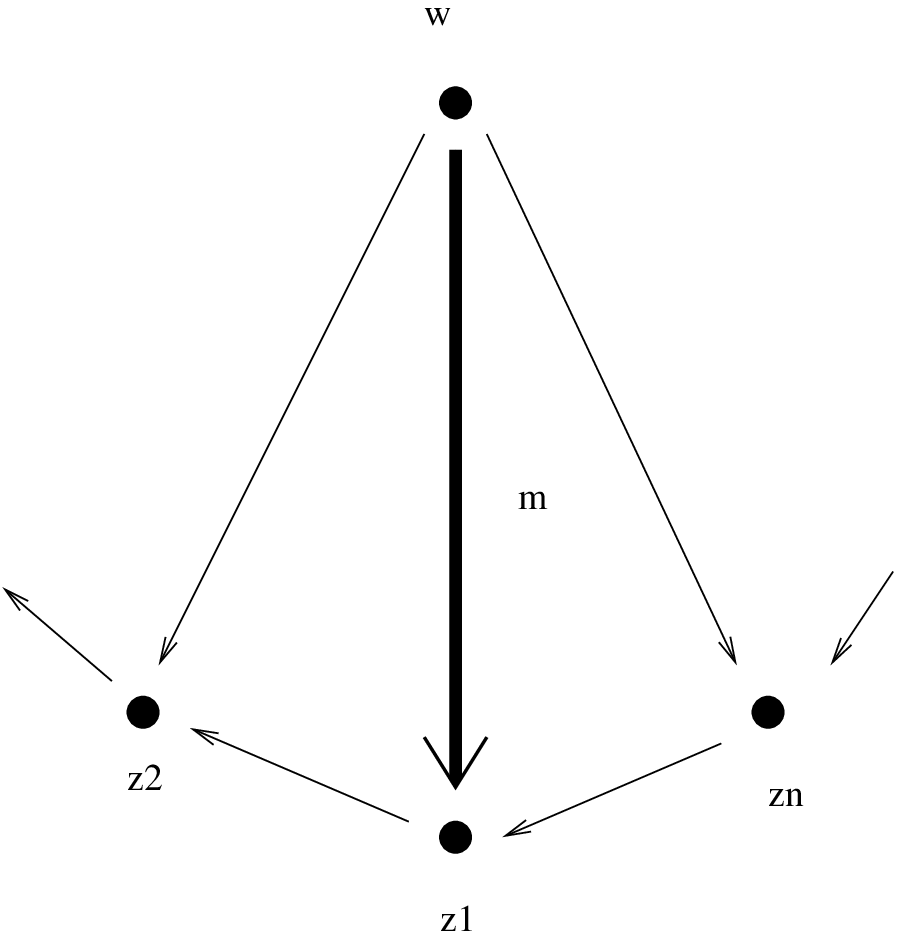}}
\pm \int \raisebox{-2.5cm}{\includegraphics[width=4cm]{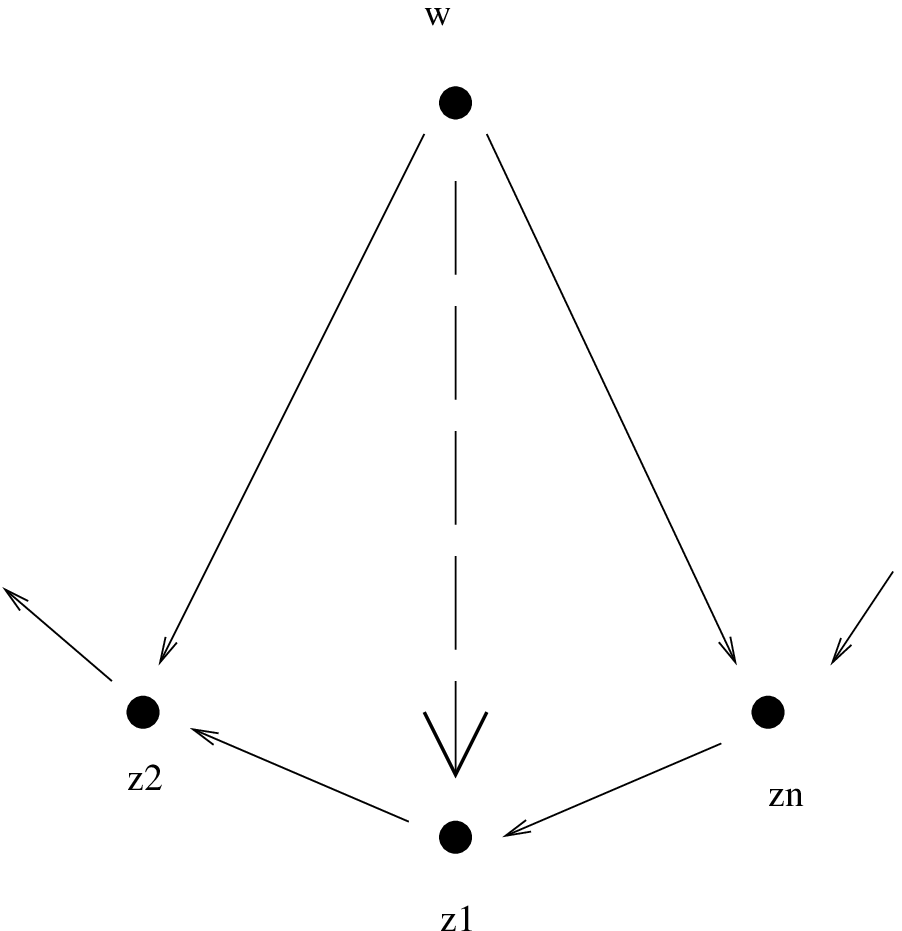}}
\end{equation}
We now simplify the integral over  $\partial_{\text{classical}}$ by looking
at the different components of $\partial_{\text{classical}}$. 
\begin{enumerate}
\item A group of vertices moves to the real line. Since every vertex
  has a least one normal outgoing arrow this means that to obtain a
  non-zero result all vertices should move to the real line which is
  excluded. So components of this type contribute nothing. 
\item A group of vertices $S$ comes together. Assume first $w\in S$.
  If $z_i\in S$, $i<n$ then to avoid double arrows after contraction
  we should also have $z_{i+1}\in S$.  If $z_1\in S$ this means we
  have to contract all vertices. I.e.\ the integral over this
  component becomes an integral over $C_{n+1}$
\[
\psfrag{z1}[][]{$z_1$}
\psfrag{z2}[][]{$z_2$}
\psfrag{zn}[][]{$z_n$}
\psfrag{m}[][]{$$}
\psfrag{w}[][]{$w$}
\int_{C_{n+1}}\raisebox{-2cm}{\includegraphics[width=4cm]{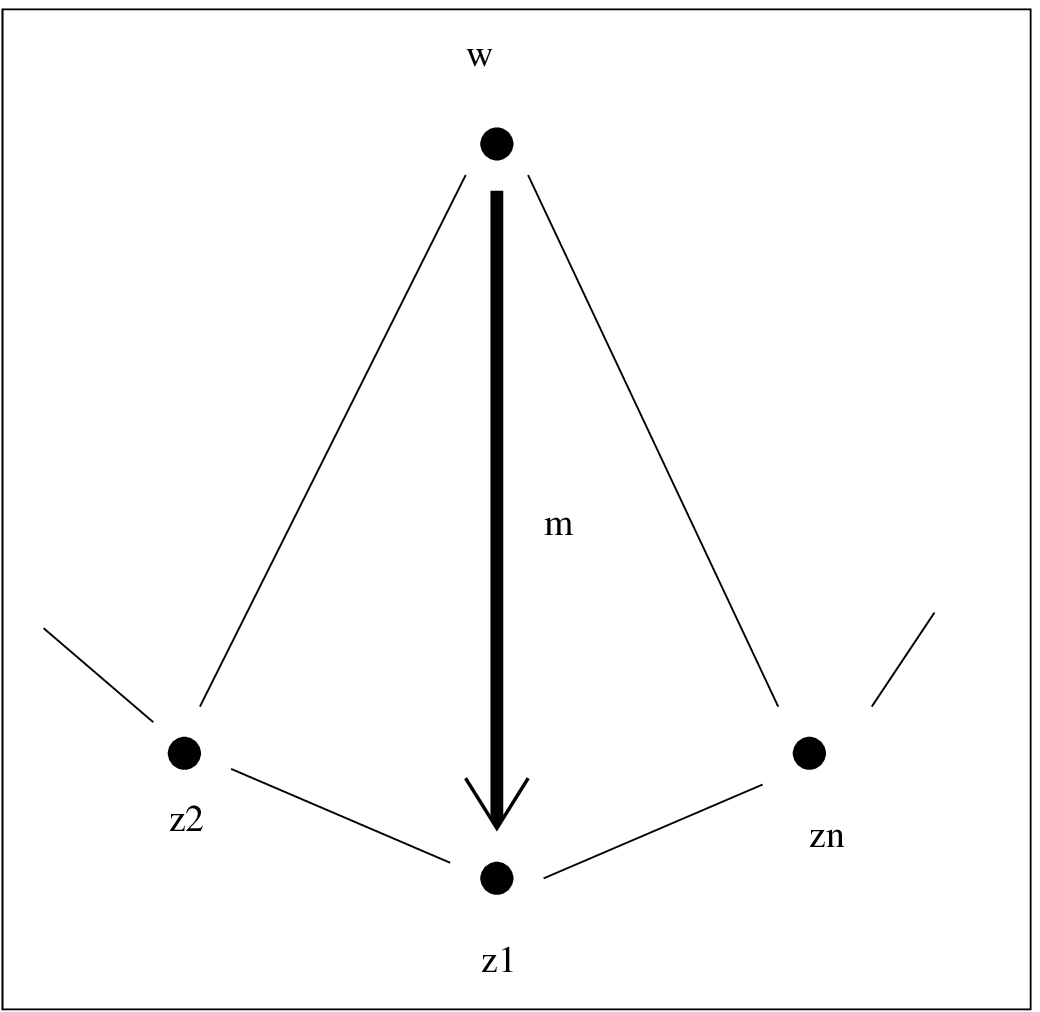}}
\]
\item
Assume $z_1\not\in S$. If $z_i\in S$ for $i<n$ then we contract more than one edge and the resulting
integral is zero by \eqref{lem66}. Thus $S=\{w,z_n\}$. The resulting 
integral is (using Rule 2)
\begin{equation}
\label{ref-3.2-7}
\psfrag{w}[][]{$w$}
\psfrag{zn}[][]{$z_n$}
\psfrag{z1}[][]{$z_1$}
\psfrag{z2}[][]{$z_2$}
\psfrag{m}[l][l]{$2$}
\psfrag{zn-1}[l][l]{$z_{n-1}$}
\frac{1}{2}\int_{C_{n,0}}
\raisebox{-2.5cm}{\includegraphics[width=4.5cm]{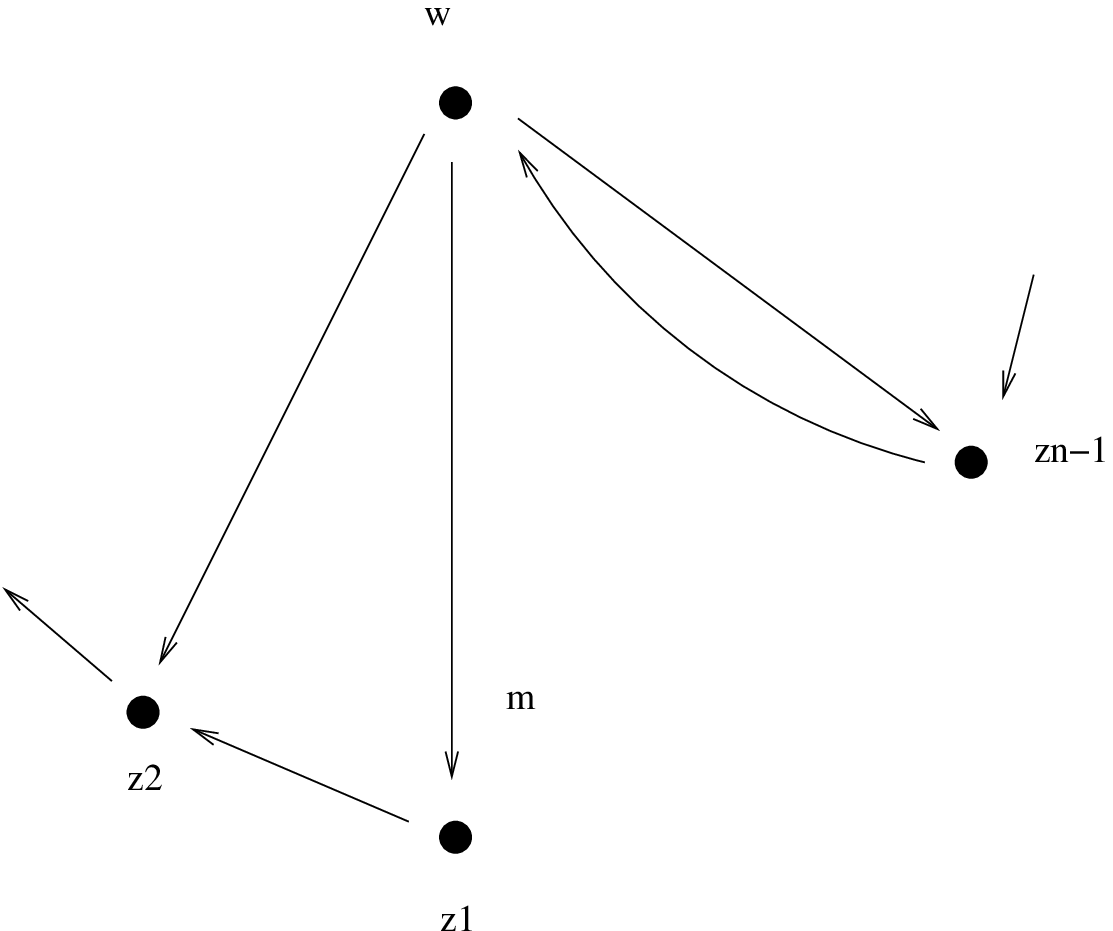}}
\end{equation}
\item Now assume $w\not \in S$. This means that we contract an edge
  $(z_i,z_{i+1})$ or $(z_n,z_1)$ (if we contract more edges then the
  integral is zero by \eqref{lem66}). The only contractions which do not
  create double edges are $(z_1,z_2)$ and~$(z_n,z_1)$. In
  that case the resulting integrals are (using Rule 2)
\begin{equation}
\psfrag{w}[][]{$w$}
\psfrag{zn}[][]{$z_n$}
\psfrag{z1}[][]{$z_1$}
\psfrag{z2}[][]{$z_2$}
\psfrag{m}[r][r]{$2$}
\psfrag{zn-1}[l][l]{$z_{n-1}$}
\frac{1}{2}\int_{C_{n.0}}
\raisebox{-2.5cm}{\includegraphics[width=4.5cm]{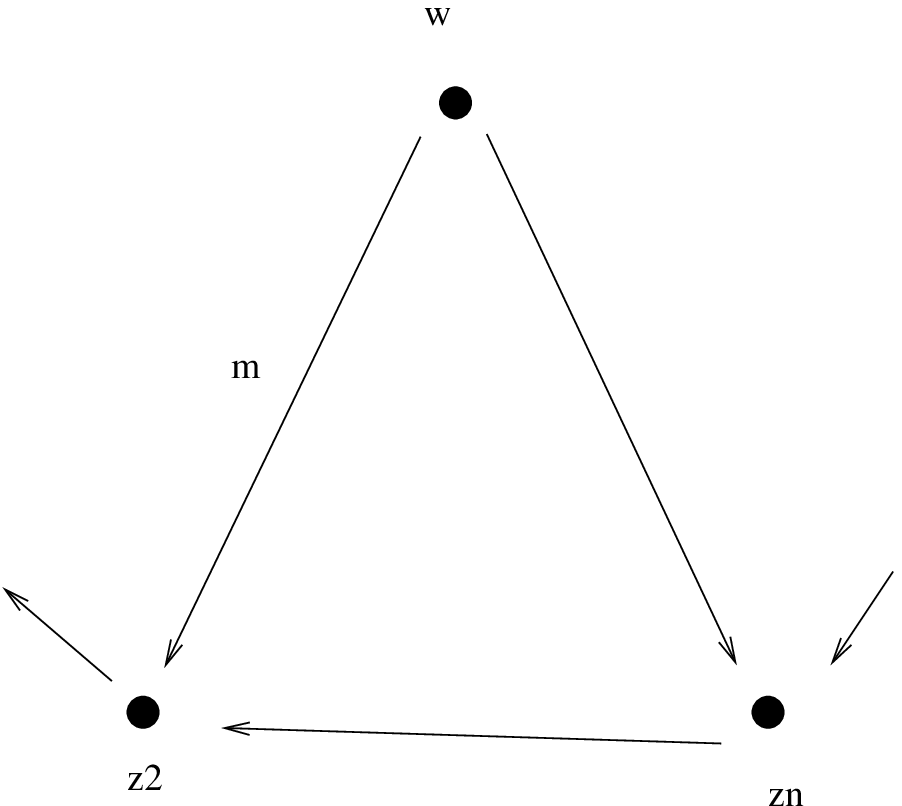}}
\end{equation}
and
\begin{equation}
\psfrag{w}[][]{$w$}
\psfrag{zn}[][]{$z_{n-1}$}
\psfrag{z2}[][]{$z_1$}
\psfrag{m}[r][r]{$2$}
\psfrag{zn-1}[l][l]{$z_{n-1}$}
\frac{1}{2}\int_{C_{n,0}}
\raisebox{-2.5cm}{\includegraphics[width=4.5cm]{wheel_bold_contract_variant.eps}}
\end{equation}
\end{enumerate}

\subsection{The boundary component associated to the branch cut}
\label{ref-3.3-8}
We need to compute the integral associated to the graph with the dashed
arrow in \eqref{ref-3.1-6}. In order to derive a recursion relation
we work more generally, starting with the following identity. 
\begin{equation}
\label{ref-3.5-9}
\psfrag{w}[][]{$w$}
\psfrag{zn}[][]{$z_n$}
\psfrag{z1}[][]{$z_1$}
\psfrag{z2}[][]{$z_2$}
\psfrag{z3}[r][r]{$z_3$}
\psfrag{m}[l][l]{$m$}
\psfrag{zn-1}[l][l]{$z_{n-1}$}
 \int \raisebox{-2.5cm}{\includegraphics[width=4cm]{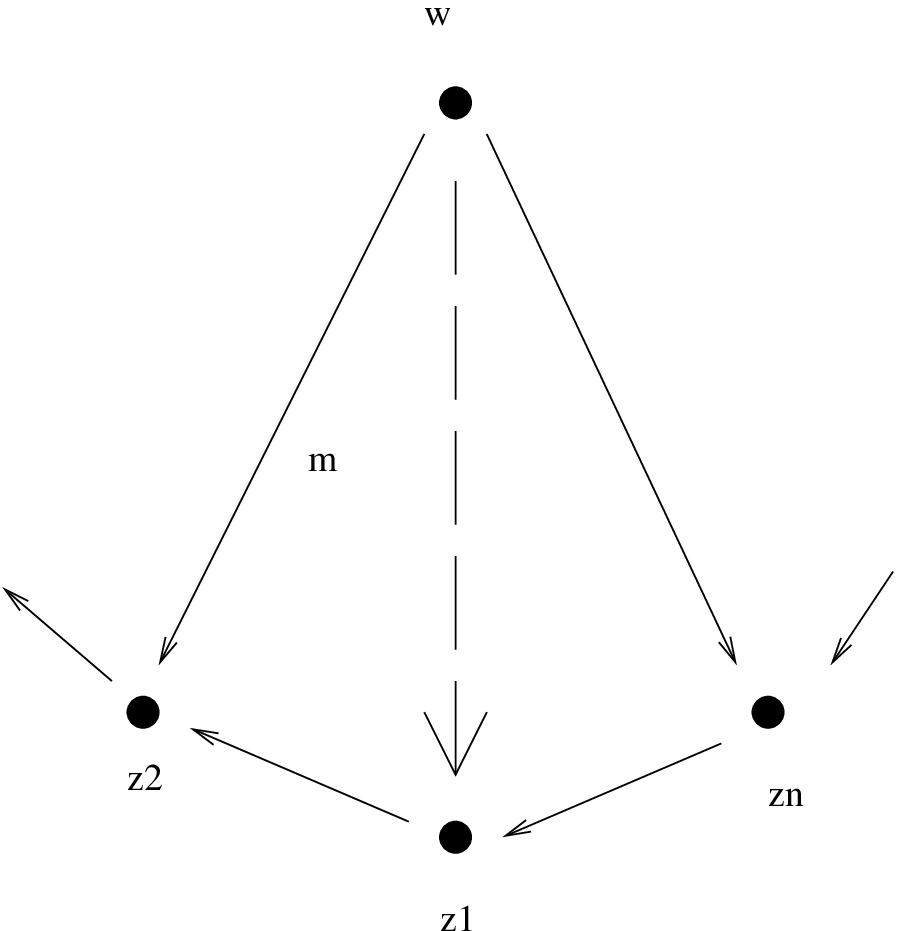}}
=\pm\int_{\partial_{\text{classical}}} \raisebox{-2.5cm}{%
\includegraphics[width=4cm]{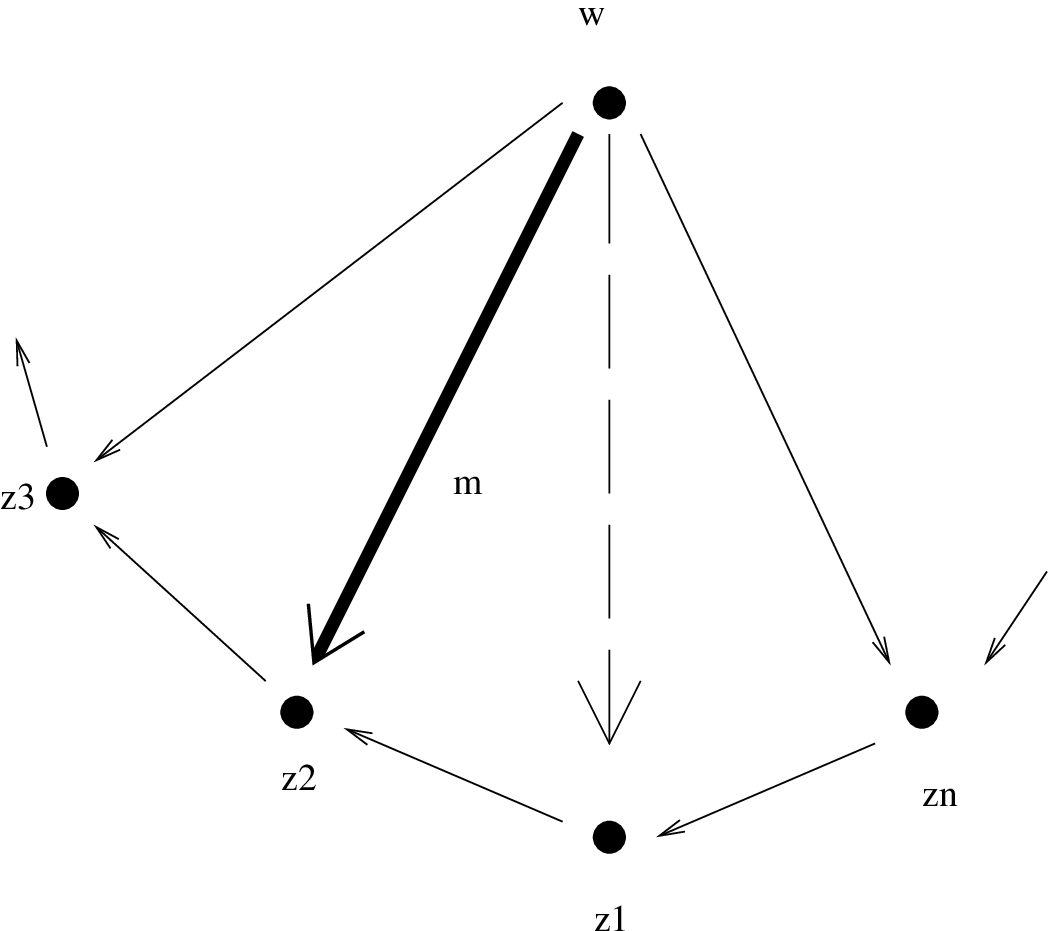}}
\pm \int \raisebox{-2.5cm}{%
\includegraphics[width=4cm]{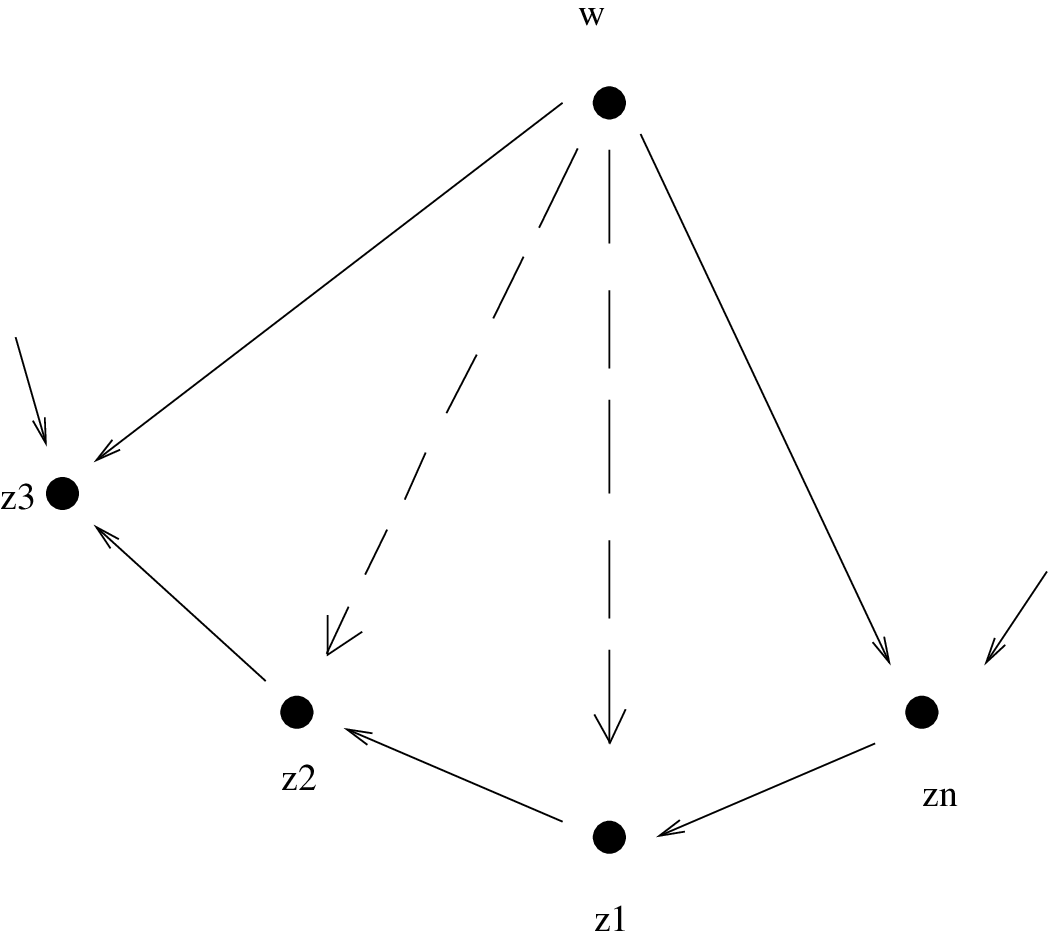}}
\end{equation}
However the integral containing the two dashed arrows is zero since  if
$z_1$, $z_2$ move on the same vertical line then $d\varphi(z_1,z_2)=0$.

It remains to determine the integral over $\partial_{\text{classical}}$. We consider
the various boundary components. We now assume $n\ge 3$.
\begin{enumerate}
\item A group of vertices moves to the real line. Since every vertex
  has at least one normal outgoing arrow (using $n\ge 2$) this means
  that to obtain a non-zero result all vertices should move to the
  real line which is excluded. So no contribution from these components.
\item A group of vertices $S$ comes together. Assume first $w\in S$.
  If $z_i\in S$, $i\neq 1,n$ then to avoid double arrows after
  contraction we should also have $z_{i+1}\in S$. However if $z_2\in
  S$ or $z_n\in S$ then $z_1\in S$ for otherwise after contraction a
  solid and a dashed arrow coincide and thus the result is zero. Thus
  we either contract all vertices or we contract
  $w,z_i,\ldots,z_n,z_1$ for $i\ge 3$ or we contract $(w,z_1)$. If we contract all vertices
  then the integral becomes
\[
\psfrag{z1}[][]{$z_1$}
\psfrag{z2}[][]{$z_2$}
\psfrag{z3}[r][r]{$z_3$}
\psfrag{zn}[][]{$z_n$}
\psfrag{m}[l][l]{$m$}
\psfrag{w}[][]{$w$}
\int \raisebox{-2cm}{\includegraphics[width=4cm]{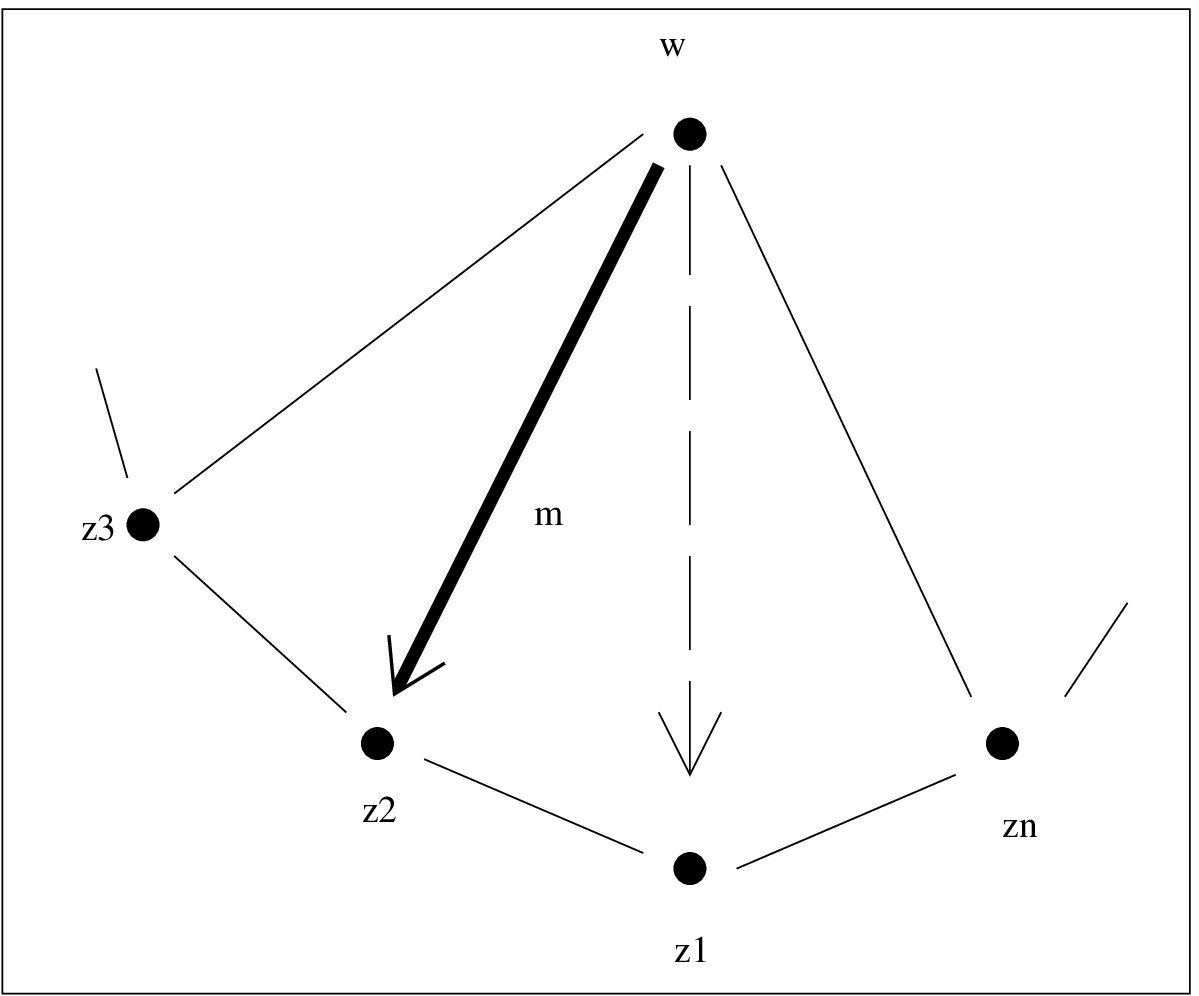}}
\]
\item Assume we contract $w,z_i,\ldots,z_n,z_1$ for $3\le i\le n$. Then the resulting
integral is \emph{a multiple of} 
\begin{equation}
\label{ref-3.6-10}
\psfrag{z1}[][]{$z_1$}
\psfrag{zi}[][]{$z_i$}
\psfrag{zi1}[l][l]{$z_{i+1}$}
\psfrag{zn}[r][r]{$z_n$}
\psfrag{m}[][]{$m$}
\psfrag{w}[][]{$w$}
\int \raisebox{-2cm}{\includegraphics[width=4cm]{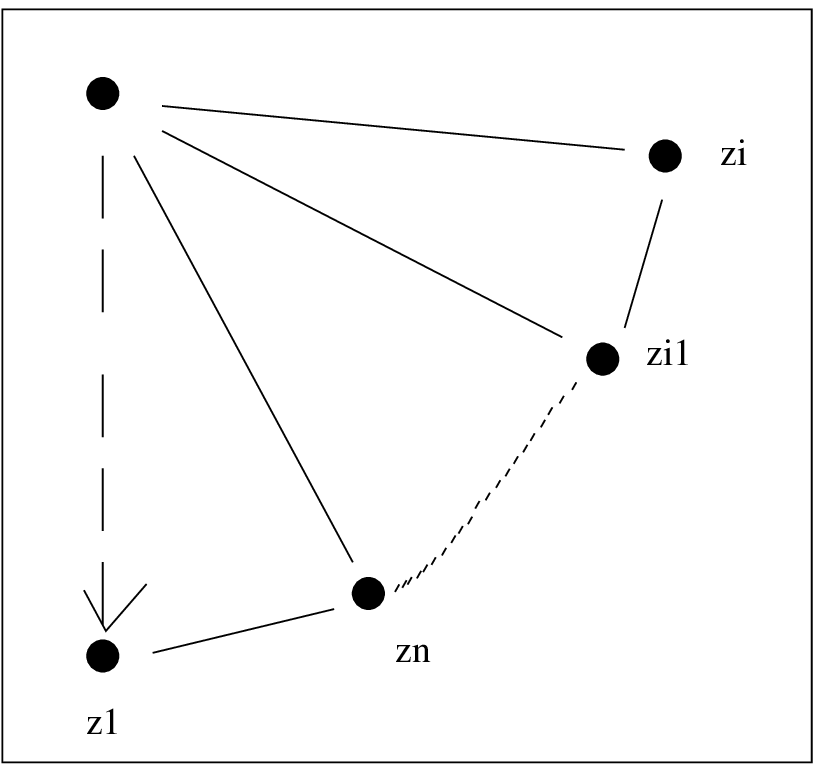}}
\end{equation}
\item It remains to consider the case $S=\{z_1,w\}$. In that case we get using 
Rule 2:
\[
\psfrag{w}[][]{$w$}
\psfrag{zn}[][]{$z_n$}
\psfrag{z1}[][]{$z_2$}
\psfrag{z2}[][]{$z_3$}
\psfrag{m}[l][l]{$m+1$}
\psfrag{zn-1}[l][l]{$z_{n}$}
\frac{1}{m+1}\int_{C_{n,0}}
\raisebox{-2.5cm}{\includegraphics[width=4.5cm]{wheel_n_bold_contracted.eps}}
\]
\item Now assume $w\not \in S$. This means that we contract an edge
  $(z_i,z_{i+1})$ or $(z_n,z_1)$ (if we contract more edges then the
  integral is zero by \eqref{lem66}). The only two possibilities
  which do not give zero immediately are $(z_1,z_2)$ and $(z_2,z_3)$. The
  resulting integrals are
\[
\psfrag{w}[][]{$w$}
\psfrag{zn}[][]{$z_n$}
\psfrag{z1}[][]{$z_2$}
\psfrag{z2}[][]{$z_3$}
\psfrag{m}[l][l]{}
\frac{1}{2} \int \raisebox{-2.5cm}{\includegraphics[width=4cm]{wheel_dashed.eps}}
\]
and (using Rule 2)
\[
\psfrag{w}[][]{$w$}
\psfrag{t}[r][r]{$m+1$}
\psfrag{zn}[][]{$z_n$}
\psfrag{z1}[][]{$z_2$}
\psfrag{z2}[][]{$z_3$}
\psfrag{m}[l][l]{}
\frac{1}{m+1} \int \raisebox{-2.5cm}{\includegraphics[width=4cm]{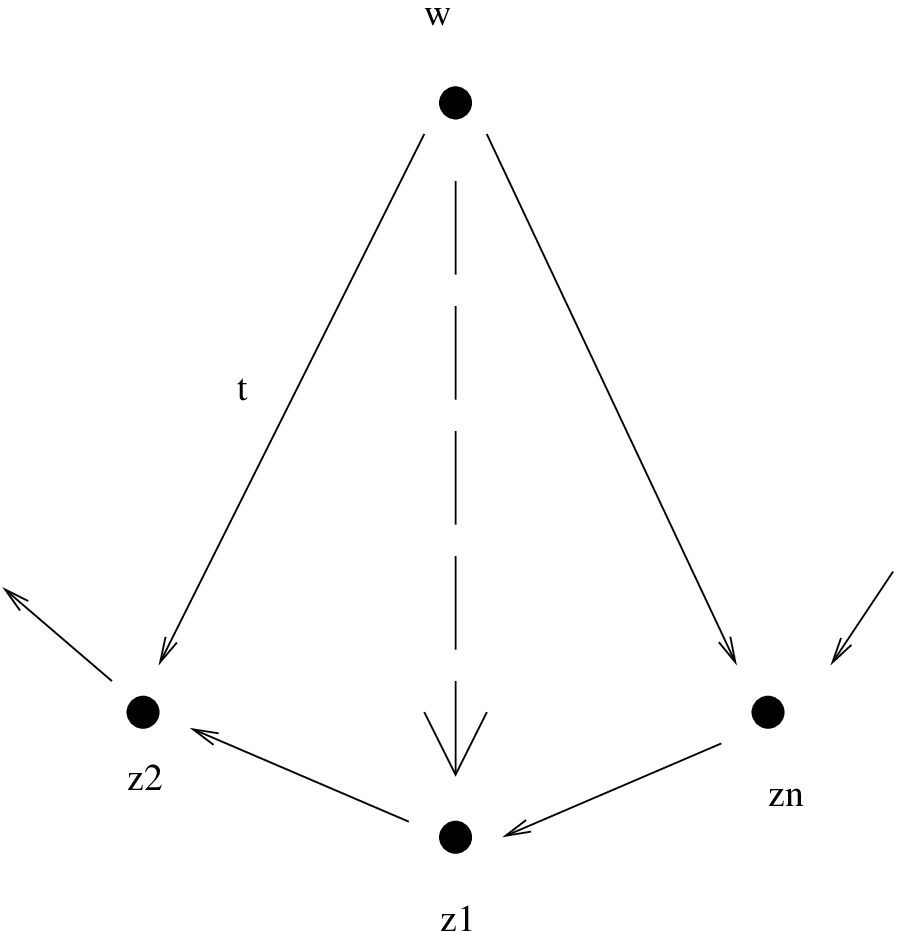}}
\]
The factor $\frac{1}{2}$ in the first integral comes from the fact
that  in this boundary component we have $\varphi(w,z_2)=0$ if $z_2$
is to the left of $z_1$.
\end{enumerate}
As will be explained below we need a second method for attacking the term with
the dashed arrow in \eqref{ref-3.1-6}. Therefore we use
\[
\psfrag{zn}[][]{$z_n$}
\psfrag{zn1}[l][l]{$z_{n-1}$}
\psfrag{z1}[][]{$z_1$}
\psfrag{z2}[][]{$z_2$}
\psfrag{m}[l][l]{$m$}
\psfrag{w}[l][l]{$w$}
\int \raisebox{-2.5cm}{\includegraphics[width=4cm]{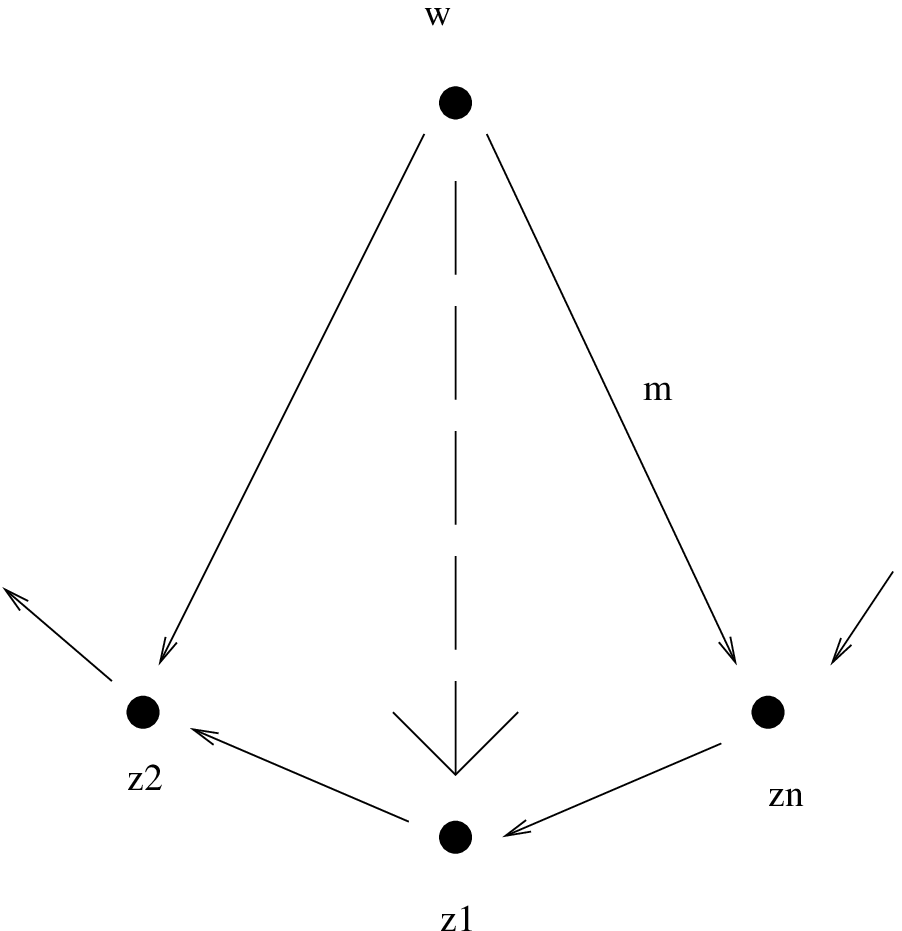}}
=\pm\int_{\partial_{\text{classical}}} \raisebox{-2.5cm}{%
\includegraphics[width=4cm]{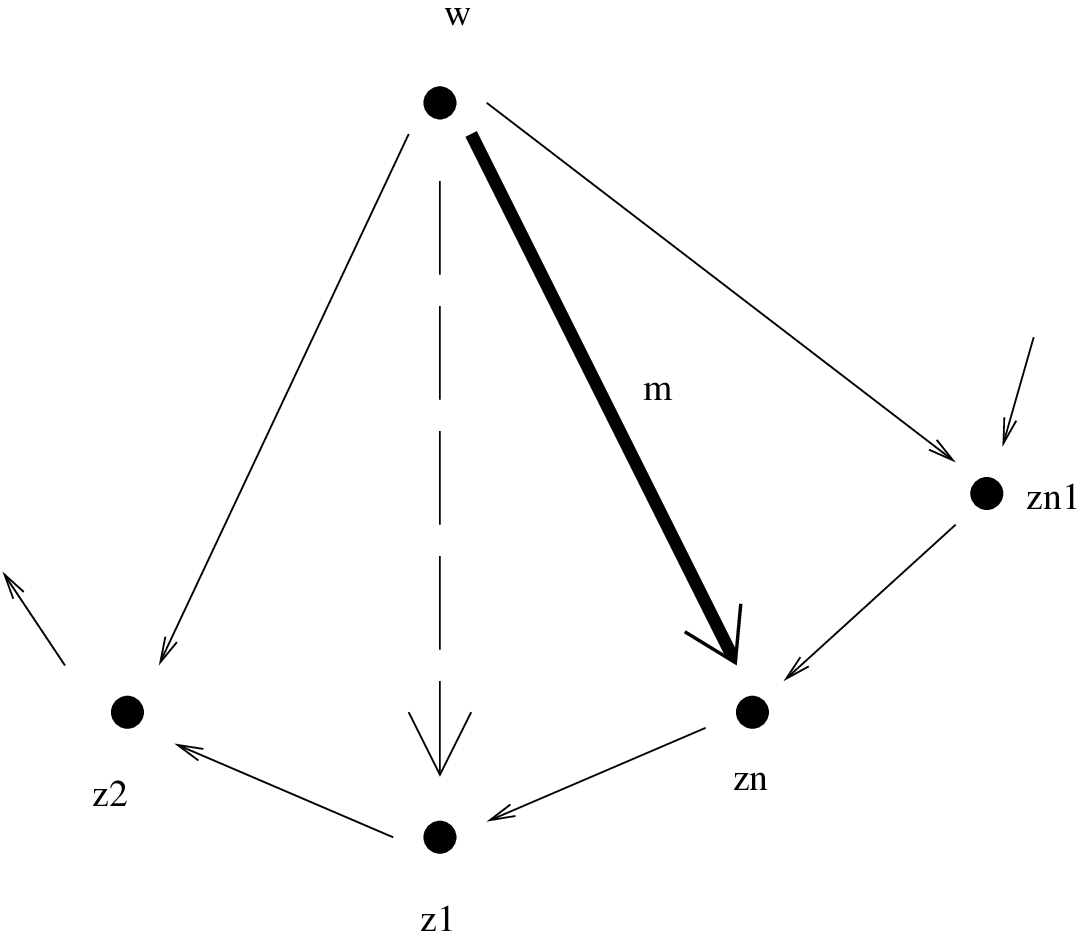}}
\pm \int \raisebox{-2.5cm}{%
\includegraphics[width=4cm]{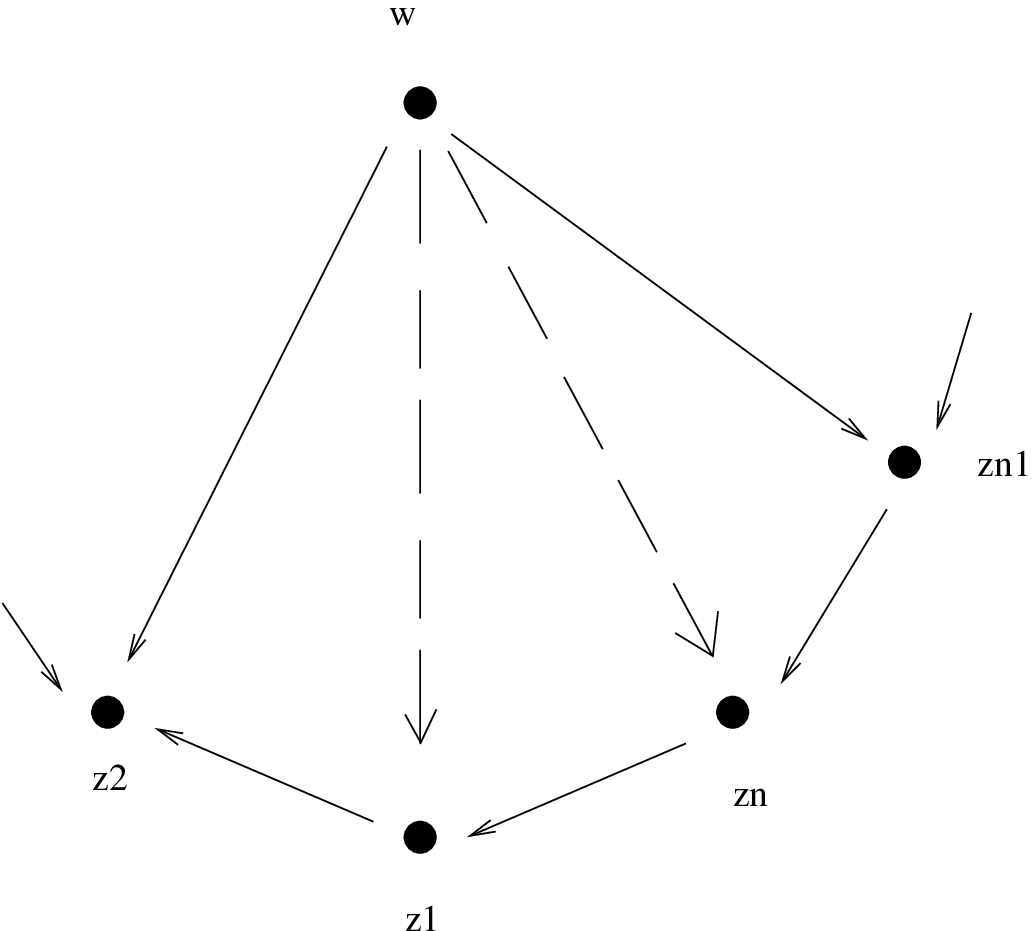}}
\]
Again the integral containing the two dashed arrows is zero. Likewise is the
integral over the boundary components where some vertices move to the real
line. We discuss the other boundary components of $\partial_{\text{classical}}$
assuming $n\ge 3$.
\begin{enumerate}
\item A group of vertices $S$ comes together. Assume first $w\in S$.
  If $z_i\in S$, $1\le i<n-1$  then to avoid double arrows after
  contraction we should also have $z_{i+1}\in S$. However if $z_n\in
  S$ of $z_2\in S$ then $z_1\in S$ for otherwise after contraction
a solid and a dashed arrow coincide and thus the result is zero. 

It follows that the possibilities for contraction are $\{z_1,\ldots,z_n,w\}$,
$\{z_1,\ldots,z_{n-1},w\}$ or
$\{z_i,\ldots,z_{n-1},w\}$ for $i\ge 3$. In  the last case we must have
$i=n-1$ by \eqref{lem66}.

If we contract all vertices then the contribution to the integral is
\[
\psfrag{zn}[][]{$z_n$}
\psfrag{zn1}[l][l]{$z_{n-1}$}
\psfrag{z1}[][]{$z_1$}
\psfrag{z2}[][]{$z_2$}
\psfrag{m}[l][l]{$m$}
\psfrag{w}[l][l]{$w$}
\int \raisebox{-1.5cm}{\includegraphics[width=4cm]{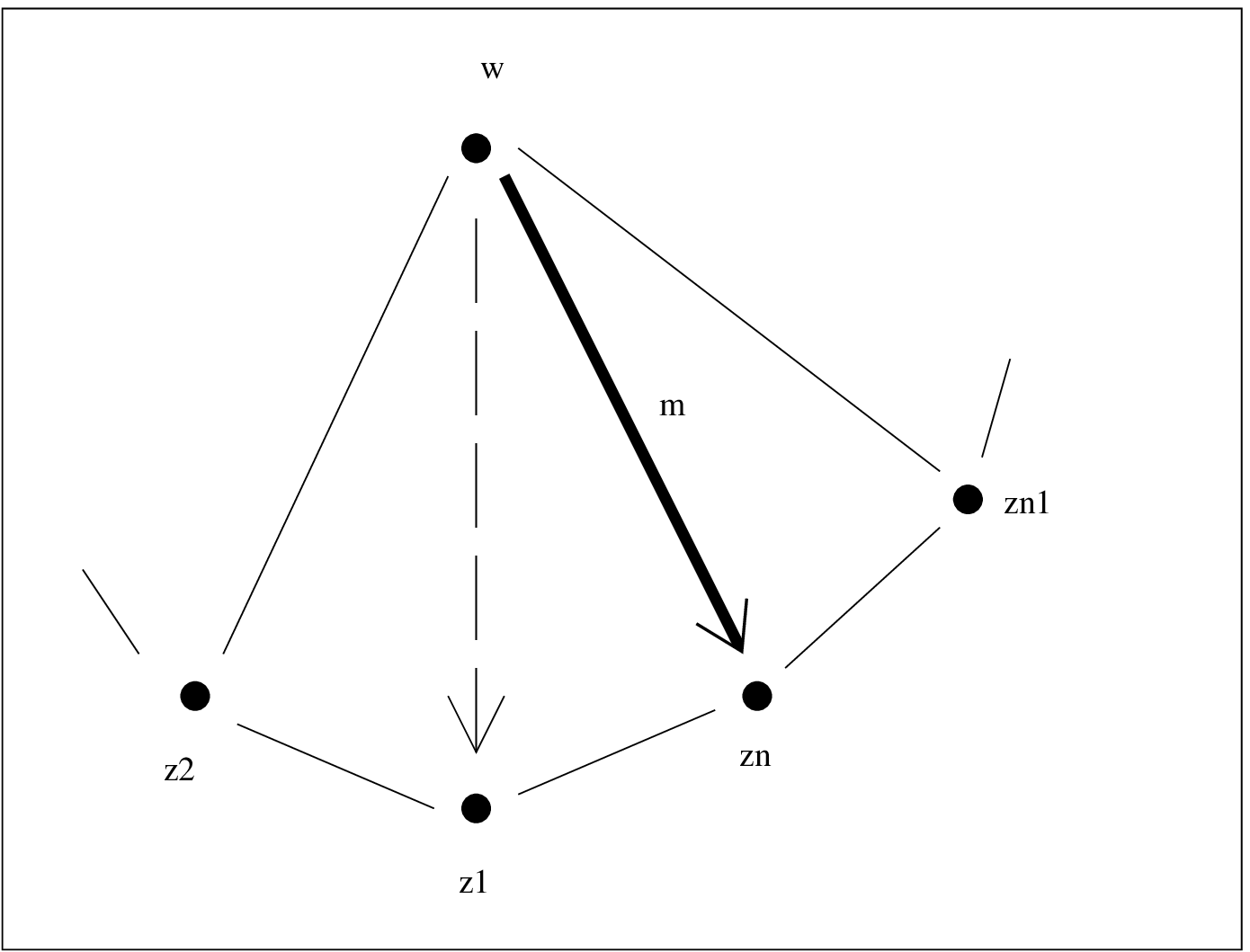}}
\]
\item If we contract
$w,z_1,\ldots,z_{n-1}$ then the resulting integral is \emph{a multiple of} 
\begin{equation}
\label{ref-3.7-11}
\psfrag{zn}[][]{$z_n$}
\psfrag{zn1}[l][l]{$z_{n-1}$}
\psfrag{zn2}[l][l]{$z_{n-2}$}
\psfrag{z1}[][]{$z_1$}
\psfrag{z2}[][]{$z_2$}
\psfrag{m}[l][l]{$m$}
\psfrag{w}[l][l]{$w$}
\int \raisebox{-1.5cm}{\includegraphics[width=4cm]{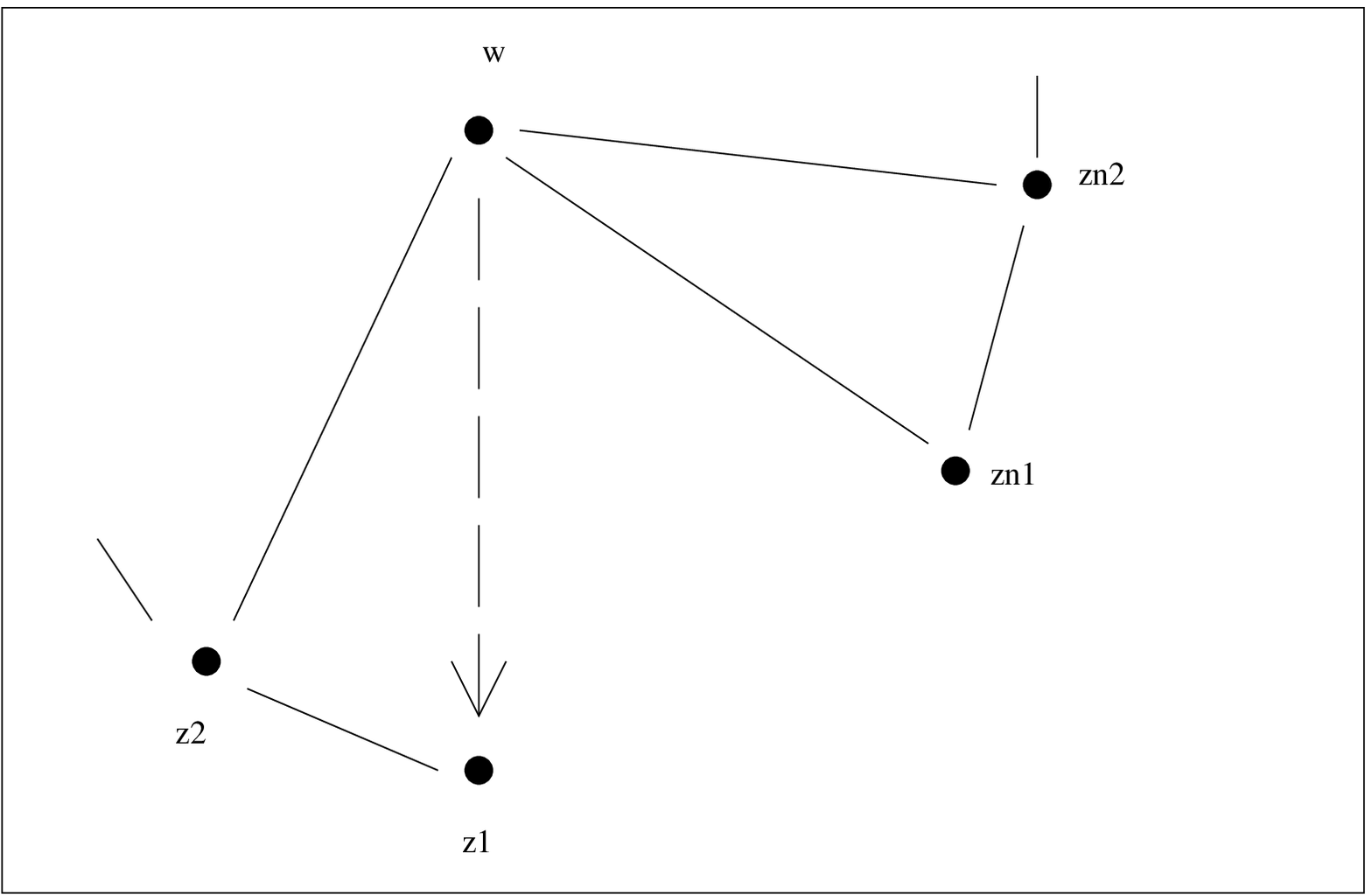}}
\end{equation}
\item It remains to consider the case
$S=\{z_{n-1},w\}$. In that case we obtain (using Rule 2)
\[
\psfrag{zn}[][]{$z_n$}
\psfrag{zn1}[l][l]{$z_{n-1}$}
\psfrag{zn2}[l][l]{$z_{n-2}$}
\psfrag{z1}[][]{$z_1$}
\psfrag{z2}[][]{$z_2$}
\psfrag{m}[l][l]{$m+1$}
\psfrag{w}[l][l]{$w$}
\frac{1}{m+1}\int \raisebox{-1.5cm}{\includegraphics[width=4cm]{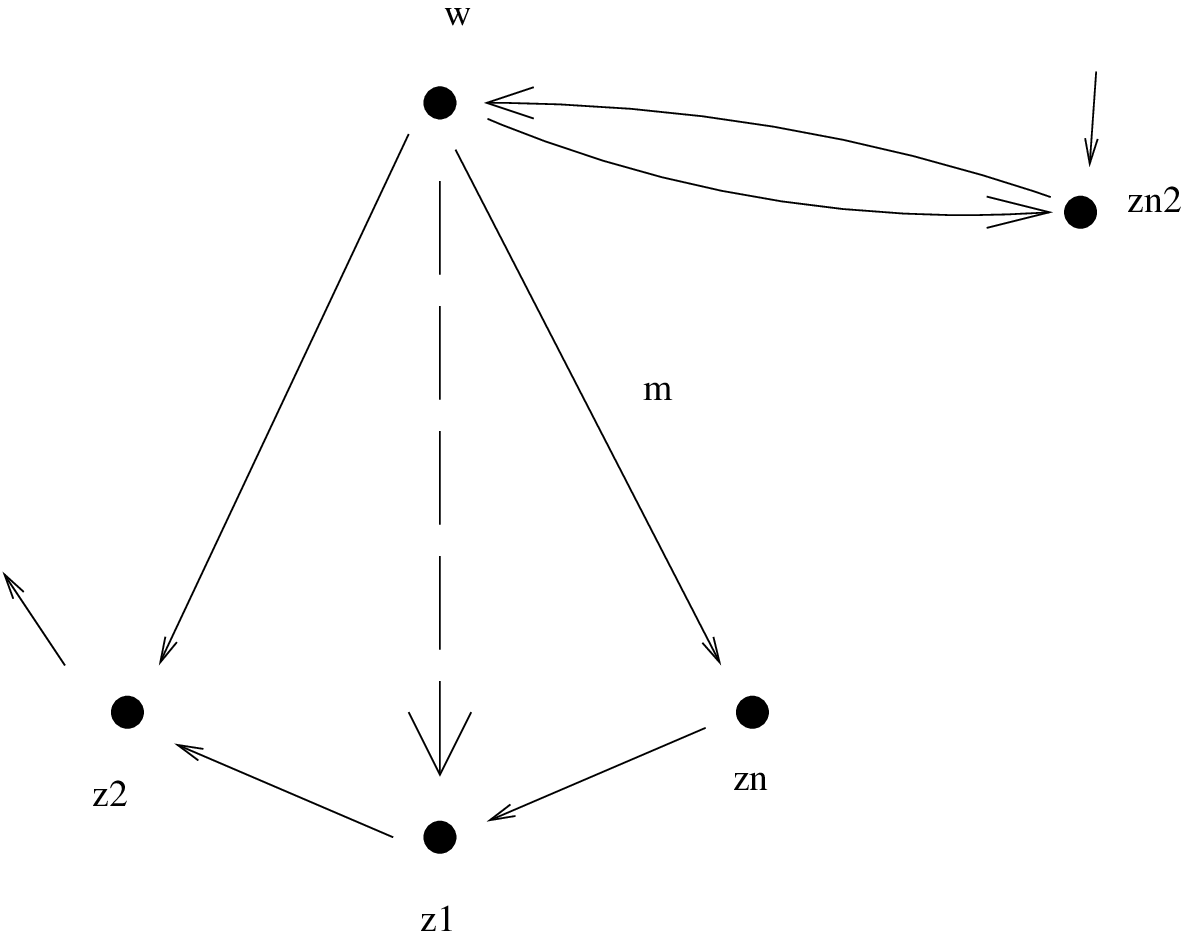}}
\]
\item 
Now assume $w\not \in S$. This means that we contract an edge
$(z_n,z_1)$ or $(z_{n-1},z_{n})$ (otherwise the result is zero by  \eqref{lem66} or by the fact that we create double arrows). The
resulting integrals are
\begin{equation}
\psfrag{w}[][]{$w$}
\psfrag{zn}[][]{$z_{n-1}$}
\psfrag{z1}[][]{$z_1$}
\psfrag{z2}[][]{$z_2$}
\psfrag{m}[r][r]{$m+1$}
\psfrag{zn-1}[l][l]{$z_{n-1}$}
\frac{1}{2}\int_{C_{n.0}}
\raisebox{-2.5cm}{\includegraphics[width=4.5cm]{wheel_dashed.eps}}
\end{equation}
and (using Rule 2)
\[
\psfrag{w}[][]{$w$}
\psfrag{m}[l][l]{$m+1$}
\psfrag{zn1}[][]{$z_{n-1}$}
\psfrag{z1}[][]{$z_2$}
\psfrag{z2}[][]{$z_3$}
\frac{1}{m+1} \int \raisebox{-2.5cm}{\includegraphics[width=4cm]{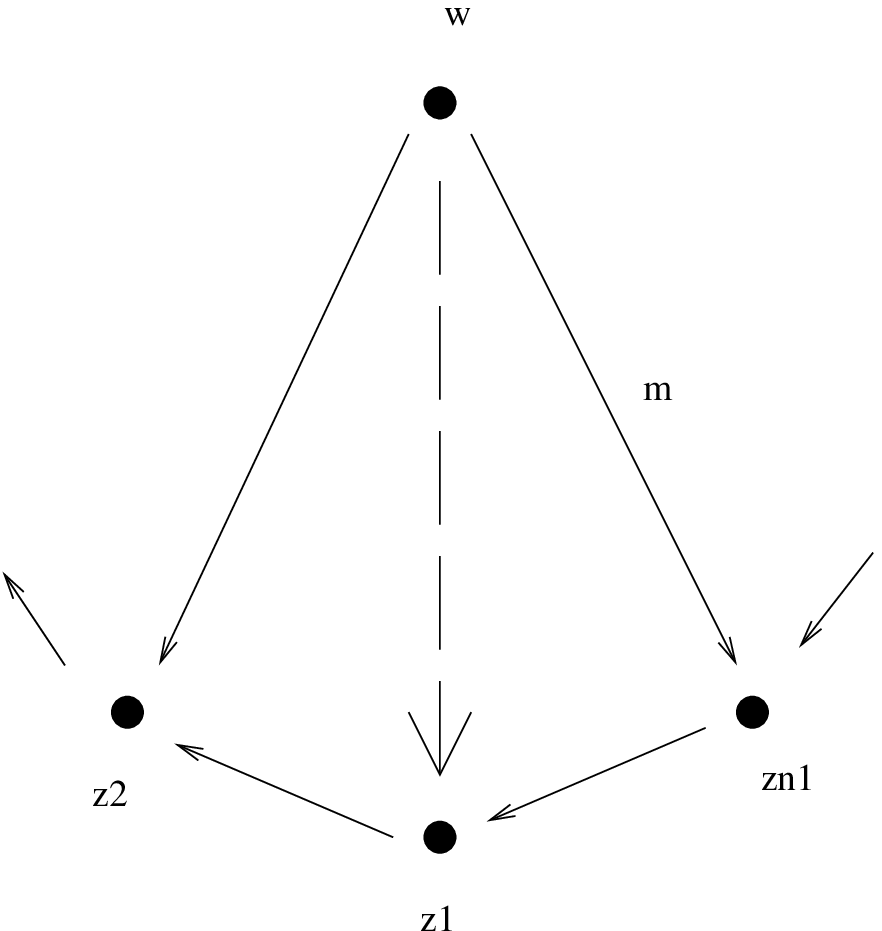}}
\]
The factor $1/2$ in the first integral comes from the fact
that due to our branch cut  $\varphi(w,z_n)=0$ if $z_2$
is to the left of $z_1$. 
\end{enumerate}

\subsection{Further comments}
\label{ref-3.4-12}
We now indicate how the relations we have derived will be used.
We still do not specify precise signs. 
\begin{enumerate}
%\item Put
%\[
%\gamma_{n,m}=\psfrag{w}[][]{$w$}
%\psfrag{zn}[][]{$z_n$}
%\psfrag{z1}[][]{$z_1$}
%\psfrag{z2}[][]{$z_2$}
%\psfrag{m}[l][l]{$m$}
%\int_{C_{n+1,0}} \raisebox{-2.5cm}{\includegraphics[width=4cm]{wheel_n.eps}}
%\]
\item We first discuss the partial evaluation of $w_n$ using
  \eqref{ref-3.1-6}. One ingredient that enters is \eqref{ref-3.2-7}. We
  define
\[
\psfrag{w}[][]{$w$}
\psfrag{zn}[][]{$z_n$}
\psfrag{z1}[][]{$z_1$}
\psfrag{z2}[][]{$z_2$}
\psfrag{m}[l][l]{$m$}
\psfrag{zn-1}[l][l]{$z_{n}$}
\alpha_{n,m}=\int_{C_{n+1,0}}
\raisebox{-2.5cm}{\includegraphics[width=4.5cm]{wheel_n_bold_contracted.eps}}
\]
We will evaluate $\alpha_{n,m}$ by subsequently integrating over
$z_1$, $z_2$, etc\dots. See \S\ref{ref-4-14} below.
\item Another ingredient that enters is \eqref{ref-8.2-51}.
\[
\psfrag{z1}[][]{$z_1$}
\psfrag{z2}[][]{$z_2$}
\psfrag{zn}[][]{$z_n$}
\psfrag{m}[][]{$$}
\psfrag{w}[][]{$w$}
\int_{C_{n+1}}\raisebox{-2cm}{\includegraphics[width=4cm]{wheel_n_bold_euclidean.eps}}
\]
We do not know how to evaluate this integral in general but we will show
by a symmetry argument that it is zero if $n$ is even (see \S\ref{ref-8.2-50}).
\item The final ingredient is given by the two identically looking integrals
\[
\psfrag{w}[][]{$w$}
\psfrag{zn}[][]{$z_n$}
\psfrag{z1}[][]{$z_1$}
\psfrag{z2}[][]{$z_2$}
\psfrag{m}[r][r]{$2$}
\psfrag{zn-1}[l][l]{$z_{n-1}$}
\int_{C_{n.0}}
\raisebox{-2.5cm}{\includegraphics[width=4.5cm]{wheel_bold_contract_variant.eps}}
\]
and
\[
\psfrag{w}[][]{$w$}
\psfrag{zn}[][]{$z_{n-1}$}
\psfrag{z2}[][]{$z_1$}
\psfrag{m}[r][r]{$2$}
\psfrag{zn-1}[l][l]{$z_{n-1}$}
\int_{C_{n,0}}
\raisebox{-2.5cm}{\includegraphics[width=4.5cm]{wheel_bold_contract_variant.eps}}
\]
We will show that the signs are such that these two integrals cancel (see \S\ref{ref-8.3-52}).
\item So ultimately we get an expression for $w_n$ in terms
  of $\alpha_{n-1,m+1}$ (which we know how to compute) and
\[
\psfrag{w}[][]{$w$}
\psfrag{zn}[][]{$z_n$}
\psfrag{z1}[][]{$z_1$}
\psfrag{z2}[][]{$z_2$}
\psfrag{m}[l][l]{$m$}
\int \raisebox{-2.5cm}{\includegraphics[width=4cm]{wheel_dashed.eps}}
\]
which is a special case of both the following integrals
\[
\beta_{n,m}\overset{\text{def}}{=}
\psfrag{w}[][]{$w$}
\psfrag{zn}[][]{$z_n$}
\psfrag{z1}[][]{$z_1$}
\psfrag{z2}[][]{$z_2$}
\psfrag{z3}[r][r]{$z_3$}
\psfrag{m}[l][l]{$m$}
\psfrag{zn-1}[l][l]{$z_{n-1}$}
\int \raisebox{-2.5cm}{\includegraphics[width=4cm]{wheel_dashed_m.eps}}
\]
and
\[
\bar{\beta}_{n,m}\overset{\text{def}}{=}
\psfrag{zn1}[l][l]{$z_{n-1}$}
\psfrag{zn}[][]{$z_n$}
\psfrag{z1}[][]{$z_1$}
\psfrag{z2}[][]{$z_2$}
\psfrag{m}[l][l]{$m$}
\psfrag{w}[l][l]{$w$}
\int \raisebox{-2.5cm}{\includegraphics[width=4cm]{wheel_dashed_m_mirror.eps}}
\]
\item We now discuss $\beta_{n,m}$ using \eqref{ref-3.5-9}. We will show
  below (see \S\ref{ref-5.2-32}) that
\[
\psfrag{z1}[][]{$z_1$}
\psfrag{zi}[][]{$z_i$}
\psfrag{zi1}[l][l]{$z_{i+1}$}
\psfrag{zn}[r][r]{$z_n$}
\psfrag{m}[][]{$m$}
\psfrag{w}[][]{$w$}
\int \raisebox{-2cm}{\includegraphics[width=4cm]{something_zero.eps}}
\]
is zero. Hence we obtain an expression for $\beta_{n,m}$ in terms of 
$\alpha_{n-1,m+1}$, $\beta_{n-1,1}$, $\beta_{n+1,m-1}$ and
\[
\epsilon_{n,m} \overset{\text{def}}{=}
\psfrag{z1}[][]{$z_1$}
\psfrag{z2}[][]{$z_2$}
\psfrag{z3}[r][r]{$z_3$}
\psfrag{zn}[][]{$z_n$}
\psfrag{m}[l][l]{$m$}
\psfrag{w}[][]{$w$}
\int \raisebox{-2cm}{\includegraphics[width=4cm]{wheel_dashed_bold_boxed.eps}}
\] 
\item In a similar way we discuss  $\bar{\beta}_{n,m}$. As already pointed out
we will show that the integral
\[
\psfrag{zn}[][]{$z_n$}
\psfrag{zn1}[l][l]{$z_{n-1}$}
\psfrag{zn2}[l][l]{$z_{n-2}$}
\psfrag{z1}[][]{$z_1$}
\psfrag{z2}[][]{$z_2$}
\psfrag{m}[l][l]{$m$}
\psfrag{w}[l][l]{$w$}
\int \raisebox{-1.5cm}{\includegraphics[width=4cm]{wheel_dashed_bold_contracted_euclidean.eps}}
\]
is zero (see \S\ref{ref-5.2-32} below).
We will also show that the integral
\[
\psfrag{zn}[][]{$z_n$}
\psfrag{zn1}[l][l]{$z_{n-1}$}
\psfrag{zn2}[l][l]{$z_{n-2}$}
\psfrag{z1}[][]{$z_1$}
\psfrag{z2}[][]{$z_2$}
\psfrag{m}[l][l]{$m+1$}
\psfrag{w}[l][l]{$w$}
\int \raisebox{-1.5cm}{\includegraphics[width=4cm]{wheel_dashed_bold_n_1_contraction.eps}}
\]
is zero (see \S\ref{ref-6.1-40}).

As a result we obtain an expression for $\bar{\beta}_{n,m}$ in terms
of  $\bar{\beta}_{n-1,m+1}$, $\beta_{n-1,1}$ and
\[
\bar{\epsilon}_{n,m}\overset{\text{def}}{=}\psfrag{zn}[][]{$z_n$}
\psfrag{zn1}[l][l]{$z_{n-1}$}
\psfrag{z1}[][]{$z_1$}
\psfrag{z2}[][]{$z_2$}
\psfrag{m}[l][l]{$m$}
\psfrag{w}[l][l]{$w$}
\int \raisebox{-1.5cm}{\includegraphics[width=4cm]{wheel_dashed_bold_m_mirror_boxed.eps}}
\]
\item We will show below by a symmetry argument that $\epsilon_{n,m}=-\bar{\epsilon}_{n,m}$ (see \S\ref{ref-7.1-44}).
\item Unfortunately we do not know how to compute $\epsilon_{n,m}$.
  However we will show below that by considering a suitable linear
  combination of $\beta_{n,m}$ and $\bar{\beta}_{n,m}$ the unknown
  quantity $\epsilon_{n,m}$ cancels out in the computation. This will ultimately
give us a recursive procedure for computing $w_n$.
\end{enumerate}
\section{A recursion relation for  $\alpha_{n,m}$}
\label{ref-4-14}
 By our definition in \S\ref{ref-3.4-12}, $\alpha_{n,m}$ is the
integral associated to the enhanced graph 
\[
\psfrag{w}[][]{$w$}
\psfrag{zn}[][]{$z_n$}
\psfrag{z1}[][]{$z_1$}
\psfrag{z2}[][]{$z_2$}
\psfrag{m}[l][l]{$m$}
\psfrag{zn-1}[l][l]{$z_{n}$}
\raisebox{-2.5cm}{\includegraphics[width=4.5cm]{wheel_n_bold_contracted.eps}}
\]
Thus with our standard ordering of edges we have
\[
\alpha_{n,m}=\frac{1}{(2\pi)^{2n+m-1}}\int_{C_{n+1,0}} d\varphi(z_1,z_2)\cdots d\varphi(z_{n-1},z_n) d\varphi(z_n,w) d\varphi(w,z_1)^m d\varphi(w,z_2)\cdots d\varphi(w,z_{n}) 
\]
In this section we prove
\begin{equation}
\label{ref-4.1-15}
\alpha_{1,m}=-\left(\frac{1}{2}-\frac{1}{m+1}\right)
\end{equation}
and for $n\ge 2$. 
\begin{equation}
\label{ref-4.2-16}
\alpha_{n,m}=(-1)^n \left( \frac{1}{2}\alpha_{n-1,2}-\frac{1}{m+1} \alpha_{n-1,m+1}\right)
\end{equation}
We put $\alpha_{0,m}=1$ such that \eqref{ref-4.2-16}
holds for $n\ge 1$. 
\subsection{Step 1} We will prove
\label{ref-4.1-17}
\begin{equation}
\label{ref-4.3-18}
\int_{z\in \Hscr\setminus\{z_1,z_2\}} 
d\varphi(z_1,z)^m d\varphi(z,z_2)=(2\pi)^m \varphi(z_1,z_2)-2\pi \varphi(z_1,z_2)^m
\end{equation}
where we have made a branch cut such that $\varphi(z_1,z)\in ]0,2\pi[$.

\begin{figure}
\psfrag{C1}[][]{$C_1$}
\psfrag{C2}[][]{$C_2$}
\psfrag{z1}[][]{$z_1$}
\psfrag{z2}[][]{$z_2$}
\psfrag{Lp}[r][r]{$B_+$}
\psfrag{Lm}[l][l]{$B_-$}
\psfrag{Cp}[][]{$C_+$}
\psfrag{Cm}[][]{$C_-$}
\psfrag{H}[][]{$H$}
\includegraphics[width=10cm]{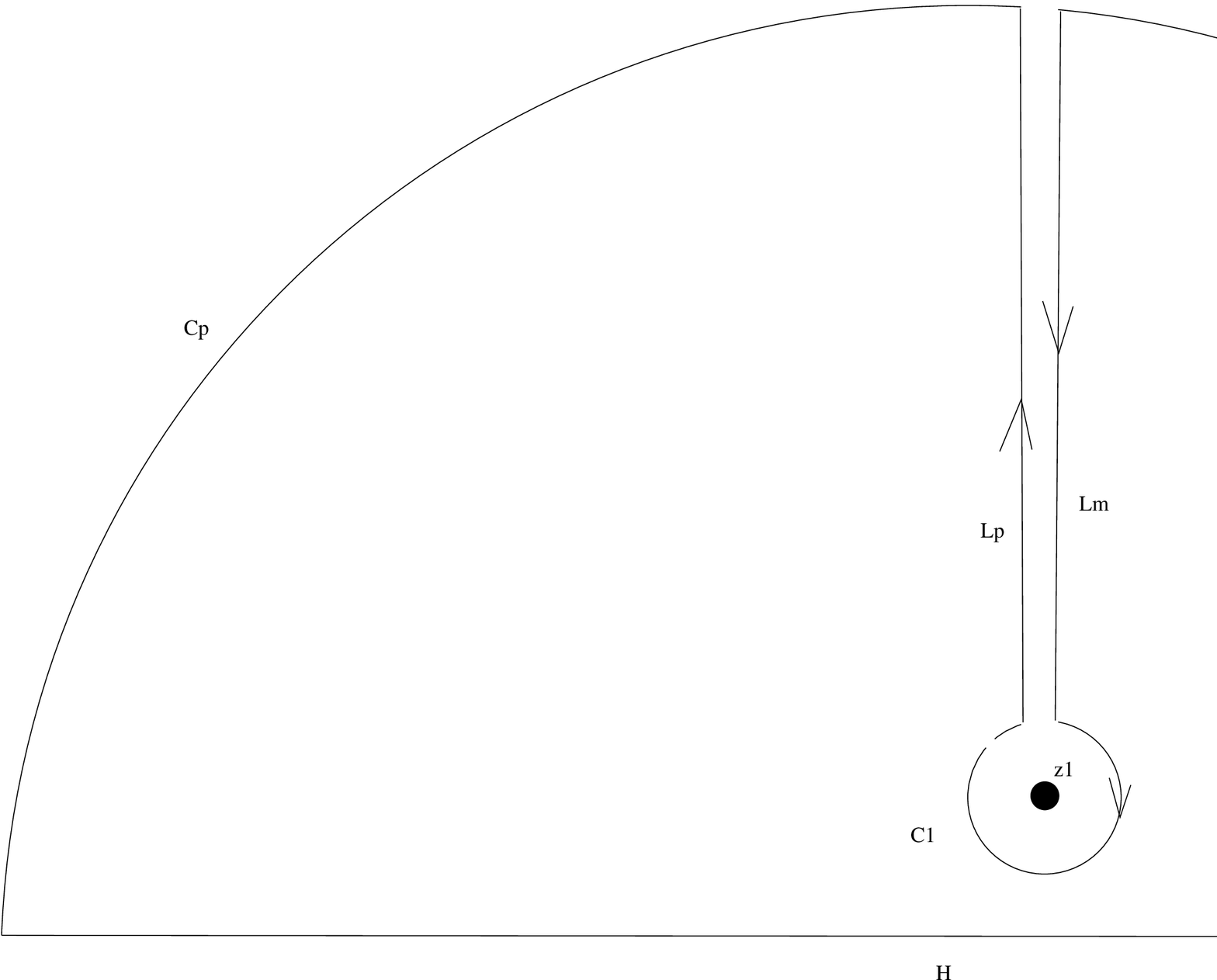}
\caption{}
\label{ref-1-13}
\end{figure}

In Figure \ref{ref-1-13}
we have indicate the integration domain and its boundary. Here $C_+\cup C_-$ is the
circle at infinity, $B_+$, $B_-$ are infinitesimally close together and $C_1$, $C_2$
have infinitesimal radius. Finally $H$ is the real line. Thus we have
\[
\int_{z\in \CC\setminus \{z_1,z_2\}} d\varphi(z_1,z)^m d\varphi(z,z_2)=
\int_{\text{boundary}} \varphi(z_1,z)^m d\varphi(z,z_2)
\]
We first compute the integral over $C_{\pm}\cup B_{\pm}\cup C_1\cup H$. 
On $B_-\cup C_-$
we have that $\varphi(z_1,z)$ is constant and equal to $2\pi$. On $B_+\cup C_+$ we
have that $\varphi(z_1,z)=0$.  On $H\cup C_1$ we have that $d\varphi(z,z_2)=0$.

Thus we find
\[
\int_{B_{\pm}\cup C_{\pm}\cup C_1\cup H}\varphi(z_1,z)^md\varphi(z,z_2)=
(2\pi)^m \int_{B_{-}\cup C_{-}}d\varphi(z,z_2)=(2\pi)^m \varphi(z_1,z_2)
\]
The integral over $C_2$ is equal to $-2\pi \varphi(z_1,z_2)^m$. This
finishes the proof of \eqref{ref-4.3-18}.
\subsection{Step 2}  
\begin{figure}
\psfrag{z1}[][]{$w$}
\psfrag{C1}[][]{$C_1$}
\psfrag{Lp}[r][r]{$B_+$}
\psfrag{Lm}[l][l]{$B_-$}
\psfrag{Cp}[][]{$C_+$}
\psfrag{Cm}[][]{$C_-$}
\psfrag{H}[][]{$H$}
\includegraphics[width=10cm]{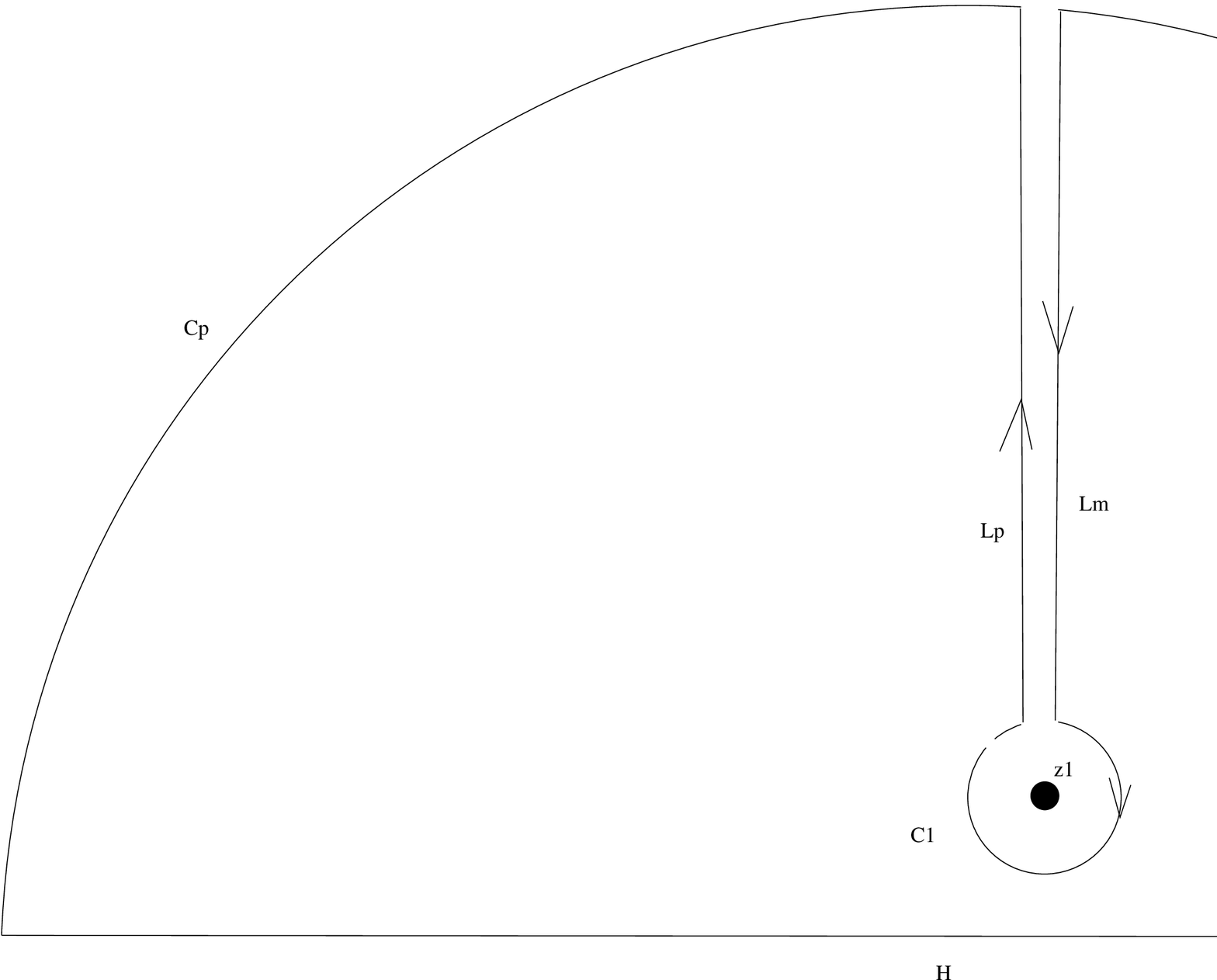}
\caption{}
\label{ref-2-19}
\end{figure}

We will prove
\label{ref-4.2-20}
\begin{equation}
\label{ref-4.4-21}
\int_{z\in \Hscr\setminus \{w\}} d\varphi(w,z)^m d\varphi(z,w)=
(2\pi)^{m+1} \left(\frac{1}{2}-\frac{1}{m+1}\right)
\end{equation}
 In Figure \ref{ref-2-19} we have indicated the integration domain and its boundary.
 Thus we
have
\[
\int_{z\in \CC\setminus \{w\}} d\varphi(w,z)^m d\varphi(z,w)=
\int_{\text{boundary}} \varphi(w,z)^m d\varphi(z,w)
\]
On $B_-\cup C_-$
we have that $\varphi(w,z)$ is constant and equal to $2\pi$. On $B_+\cup C_+$ we
have that $\varphi(w,z)=0$.  On $H$ we have that $d\varphi(z,w)=0$.
Thus
\[
\int_{C_{\pm}\cup B_{\pm}\cup H} \varphi(w,z)^m d\varphi(z,w)=(2\pi)^m\pi
\]
It remains to compute
\[
\int_{C_1} \varphi(w,z)^m d\varphi(z,w)
\]
We put $z=w+re^{-i\theta}$ for $\theta\in ]-\pi/2,3\pi/2[$. We have
$\varphi(w,z)=3\pi/2-\theta$ and hence $d\varphi(z,w)=-d\theta$. Thus
we get
\[
\int_{C_1} \varphi(w,z)^m d\varphi(z,w)=-\int_{-\pi/2}^{3\pi/2}  (3\pi/2-\theta)^md\theta
=-\int_0^{2\pi}\alpha^m d\alpha=-\frac{(2\pi)^{m+1}}{m+1}
\] 
This finishes the proof of \eqref{ref-4.4-21}.
\subsection{Step 3} Now we prove \eqref{ref-4.1-15}.
\label{ref-4.3-22}
\[
\alpha_{1,m}=\frac{1}{(2\pi)^{m+1}}\int d\varphi(z_1,w) d\varphi(w,z_1)^m=
-\frac{1}{(2\pi)^{m+1}}\int  d\varphi(w,z_1)^m d\varphi(z_1,w)=
-\left(\frac{1}{2}-\frac{1}{m+1}\right)
\]
by \eqref{ref-4.4-21}. 
\subsection{Step 4}
\label{ref-4.4-23}
Now we prove \eqref{ref-4.2-16}.
\[
\alpha_{n,m}=(-1)^{n}\frac{1}{(2\pi)^{2n+m-1}}\int d\varphi(w,z_1)^m d\varphi(z_1,z_2) \cdots d\varphi(z_{n-1},z_n) d\varphi(z_n,w) d\varphi(w,z_2)\cdots d\varphi(w,z_{n})
\]
Integrating over $z_1$ and using \eqref{ref-4.3-18} first we find
\begin{multline*}
\alpha_{n,m}=(-1)^{n}\frac{1}{(2\pi)^{2n+m-1}}\int ((2\pi)^m \varphi(w,z_2)-
2\pi \varphi(w,z_2)^m)\times\\
d\varphi(z_2,z_3)\cdots d\varphi(z_{n-1},z_n) d\varphi(z_n,w) d\varphi(w,z_2)\cdots d\varphi(w,z_{n})
\end{multline*}
which is equal to
\begin{multline*}
(-1)^{n}\frac{1}{2(2\pi)^{2n-1}}\int  d\varphi(z_2,z_3)\cdots d\varphi(z_{n-1},z_n) d\varphi(z_n,w)  d\varphi(w,z_2)^2   \cdots d\varphi(w,z_{n-1})\\
-
(-1)^{n}\frac{1}{(m+1)(2\pi)^{2n+m-2}} \int d\varphi(z_2,z_3)\cdots d\varphi(z_{n-1},z_n) d\varphi(z_n,w) d\varphi(w,z_2)^{m+1}\cdots d\varphi(w,z_{n})
\end{multline*}
So that we indeed find
\[
\alpha_{n,m}=(-1)^n \left( \frac{1}{2}\alpha_{n-1,2}-\frac{1}{m+1} \alpha_{n-1,m+1}\right)
\]
\section{A recursion relation for $\beta_{n,m}$}
\label{ref-5-24}
 By our definition in \S\ref{ref-3.4-12}, $\beta_{n,m}$ is the
integral associated to the enhanced graph 
\[
\psfrag{z1}[][]{$z_1$}
\psfrag{z2}[][]{$z_2$}
\psfrag{zn}[][]{$z_n$}
\psfrag{m}[l][l]{$m$}
\psfrag{w}[][]{$w$}
\includegraphics[width=4cm]{wheel_dashed_m.eps}
\]
Thus with our standard edge ordering
\[
\beta_{n,m}=
\frac{1}{(2\pi)^{2n+m-2} }
\int_{B_-\setminus D}
 d\varphi(z_1,z_2)\cdots d\varphi(z_{n-1},z_n)d\varphi(z_n,z_{1})
 d\varphi(w,z_2)^m\cdots d\varphi(w,z_{n}) 
\]
where
\[
B_-=\{\Re z_1=\Re w+, \Im {z_1}>\Im w\}\subset C_{n+1,0}
\]
\[
D=\{\Re z_2=\Re w, \Im {z_2}>\Im w\}\subset C_{n+1,0}
\]
($D$ represents the branch cut associated to the arrow
$w\xrightarrow{m} z_2$, see \S\ref{ref-3.1-4}).  We need to be careful
about specifying the orientation of the integration domain. We do this
next. We first normalize things by putting $w=i$.  If $z_j=x_j+iy_j$
then $C_{n+1,0}$ is oriented by \cite{AMM1}
\[
dx_1dy_1dx_2dy_2\cdots dx_ndy_n
\]
In the neighborhood of $B_-$ we have coordinates
$(x_1,y_1,\ldots,x_n,y_n)$ with $x_1\ge 0$. Therefore if $B_-$ is
oriented with its outgoing normal (as is necessary for the application
of Stokes theorem) then it is oriented by
\begin{equation}
\label{ref-5.1-25}
-dy_1dx_2dy_2\cdots dx_ndy_n
\end{equation}
Having determined the orientation of $B_-$ we let $w$ be free again. In this section we will prove 
\begin{equation}
\label{ref-5.2-26}
\beta_{2,m}=-\frac{1}{8}+\frac{1}{m+1}\left( \frac{1}{2}-\frac{1}{m+2}\right)
\end{equation} 
and for $n\ge 3$ 
\begin{equation}
\label{ref-5.3-27}
\beta_{n,m}=-\frac{1}{m+1}\alpha_{n-1,m+1}
+(-1)^{n+1} \frac{1}{2}\beta_{n-1,1}+(-1)^{n} \frac{1}{m+1} \beta_{n-1,m+1}
+(-1)^{n}\epsilon_{n,m}
\end{equation}
where $\epsilon_{n,m}$ is the integral associated to the graph
\[
\psfrag{z1}[][]{$z_1$}
\psfrag{z2}[][]{$z_2$}
\psfrag{z3}[r][r]{$z_3$}
\psfrag{zn}[][]{$z_n$}
\psfrag{m}[l][l]{$m$}
\psfrag{w}[][]{$w$}\raisebox{-2cm}{\includegraphics[width=4cm]{wheel_dashed_bold_boxed.eps}}
\]
(see \S\ref{ref-3.4-12}). Thus with our standard orientation on edges
\[
\epsilon_{n,m}=
\frac{1}{(2\pi)^{2n+m-2}}\int \theta(w,z_2)^m d\theta(z_1,z_2)d\theta(z_2,z_3)
\cdots d\theta(z_n,z_1) d\theta(w,z_3)\cdots d\theta(w,z_n)
\]
Here we put $w=0$, $z_1=i$ and the integral is over the complement
of the diagonal in $(\CC-\{0,i\})^{n-1}$.
\subsection{Step 1} 
\begin{figure}
\psfrag{L1}[r][r]{$L_1$}
\psfrag{L2}[r][r]{$L_2$}
\psfrag{C1}[l][l]{$C_1$}
\psfrag{C}[l][l]{$C$}
\psfrag{z1}[r][r]{$z_1$}
\psfrag{H}[r][r]{$H$}
\includegraphics[width=8cm]{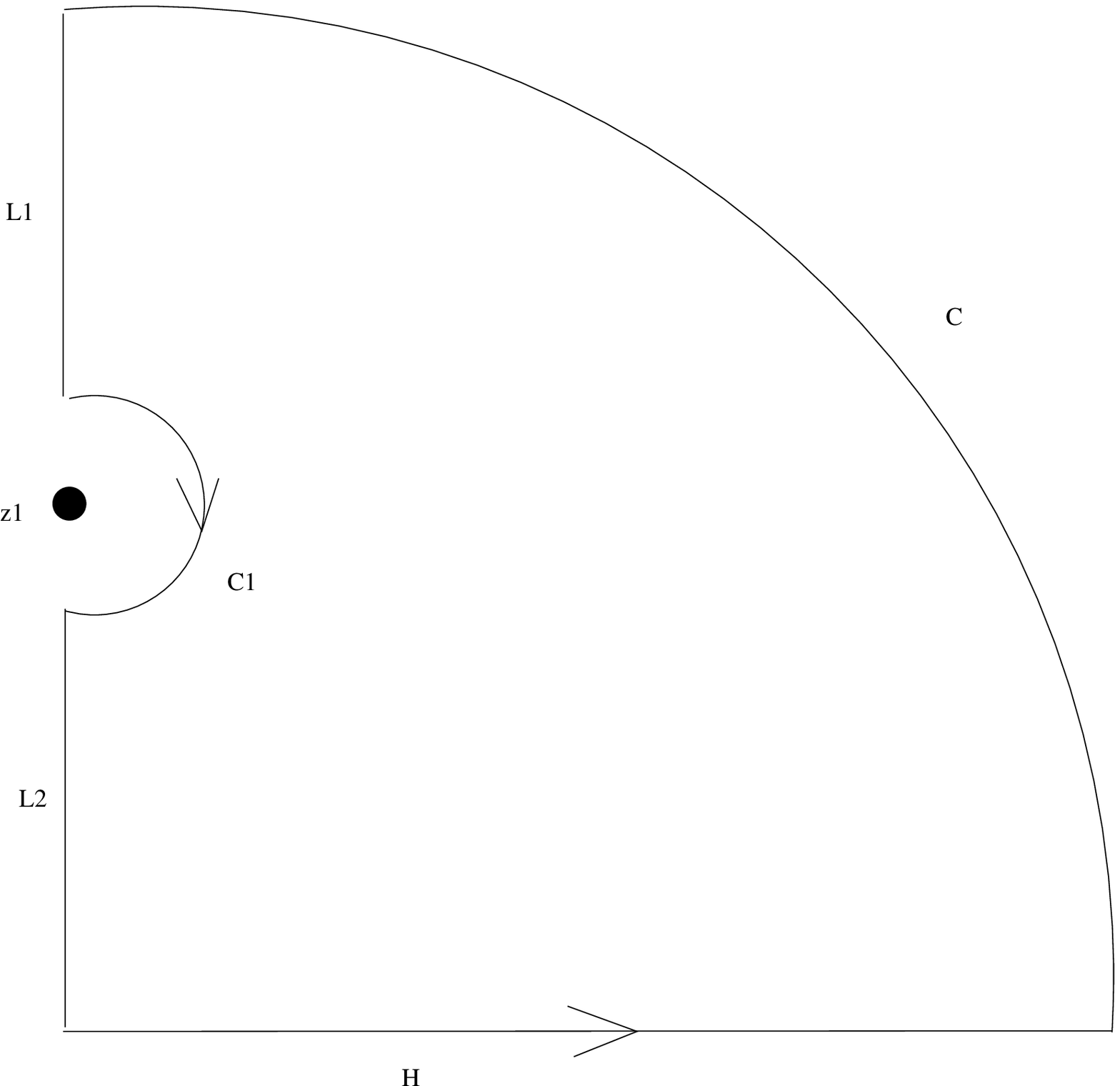}
\caption{}
\label{ref-3-30}
\end{figure}

We first consider the case $n=2$. Thus
\label{ref-5.1-28}
\[
\beta_{2,m}=
\frac{1}{(2\pi)^{2+m}}
\int_{B_-\setminus D}
 d\varphi(z_1,z_2)d\varphi(z_{2},z_1)
 d\varphi(w,z_2)^m
\]
We first integrate over $w$. We get
\[
\int_{B_-\setminus D}  d\varphi(w,z_2)^m=
\varphi(z_1,z_2)^m-(2\pi)^m [z_1,z_2]
\]
where
\[
[z_1,z_2]=
\begin{cases}
1&\text{$z_1$ is to the left of $z_2$}\\
0&\text{otherwise}
\end{cases}
\]
Thus we must compute
\begin{equation}
\label{ref-5.4-29}
\frac{1}{(2\pi)^{2+m}}
\int 
\varphi(z_1,z_2)^m d\varphi(z_1,z_2)d\varphi(z_{2},z_1)
\end{equation}
and 
\begin{equation}
\label{ref-5.5-31}
\frac{1}{(2\pi)^{2}}
\int [z_1,z_2] d\varphi(z_1,z_2)d\varphi(z_{2},z_1)
\end{equation}
Using Stokes theorem we rewrite \eqref{ref-5.5-31} as a path integral as
in Figure \ref{ref-3-30}
\[
\frac{1}{(2\pi)^{2}}
\int_{L_1\cup C_1\cup L_2\cup H\cup C} \varphi(z_1,z)d\varphi(z,z_1)
\]
If $z\in L_{1,2}\cup H$ then $d\varphi(z,z_1)=0$. If $z\in C$ then
$\varphi(z_1,z)=2\pi$. 

Thus
\[
\frac{1}{(2\pi)^{2}} \int_{L_1\cup L_2\cup H\cup C}
\varphi(z_1,z)d\varphi(z,z_1)=\frac{1}{(2\pi)}(\pi-0)=\frac{1}{2}
\]
Now we compute
\[
\frac{1}{(2\pi)^{2}} \int_{C_1}
\varphi(z_1,z)d\varphi(z,z_1)
\]
We put $z=z_1+re^{-i\theta}$ for $\theta\in ]-\frac{\pi}{2},\frac{\pi}{2}[$. 
Then $\varphi(z_1,z)=\frac{3}{2}\pi-\theta$  and hence $d\varphi(z,z_1)=-d\theta$.

Thus
\begin{align*}
\frac{1}{(2\pi)^{2}} \int_{C_1}
\varphi(z_1,z)d\varphi(z,z_1)
&=
-\frac{1}{(2\pi)^{2}}\int_{-\pi/2}^{\pi/2} (\frac{3}{2}\pi-\theta)d\theta\\
&=\frac{1}{(2\pi)^{2}}\int_{-\pi/2}^{\pi/2} (\theta-\frac{3}{2}\pi)d\theta\\
&=\frac{1}{(2\pi)^{2}}\int_{-2\pi}^{-\pi} \theta d\theta\\
&=-\frac{3}{8}
\end{align*}
and hence
 \[
\frac{1}{(2\pi)^{2}}
\int [z_1,z_2] d\varphi(z_1,z_2)d\varphi(z_{2},z_1)=\frac{1}{2}-\frac{3}{8}=\frac{1}{8}
\]
One the other hand by \eqref{ref-4.3-18}
\begin{align*}
\frac{1}{(2\pi)^{2+m}}
\int 
\varphi(z_1,z_2)^m d\varphi(z_1,z_2)d\varphi(z_{2},z_1)&=
\frac{1}{(m+1)(2\pi)^{2+m}}
\int 
\ d\varphi(z_1,z_2)^{m+1}d\varphi(z_{2},z_1)\\ &=
\frac{1}{m+1}\left( \frac{1}{2}-\frac{1}{m+2}\right)
\end{align*}
so that our final formula is 
\[
\beta_{2,m}=-\frac{1}{8}+\frac{1}{m+1}\left( \frac{1}{2}-\frac{1}{m+2}\right)
\]
\subsection{Step 2} We claim
\label{ref-5.2-32}
\begin{equation}
\label{ref-5.6-33}
\int_{z\in \CC\setminus \{z_1,z_2\}} d\theta(z_1,z)d\theta(z,z_2)=0
\end{equation}
Since $d\theta(z,z_2)=d\theta(z_2,z)$ we may as well prove
\[
\int_{z\in \CC\setminus \{z_1,z_2\}} d\theta(z_1,z)d\theta(z_2,z)=0
\]
We may assume $z_1\neq z_2$ since otherwise the claim is trivial. 
 Consider the following figure
\[
\psfrag{z1}[][]{$z_1$}
\psfrag{z2}[][]{$z_2$}
\psfrag{z}[][]{$z$}
\psfrag{zp}[][]{$z'$}
\includegraphics[width=4cm]{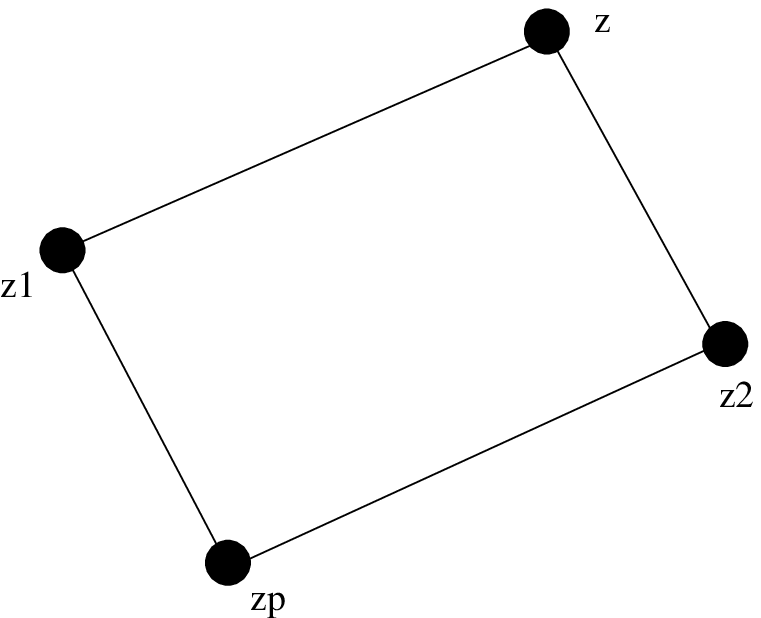}
\]
Then we have $d\theta(z_1,z)d\theta(z_2,z)=-d\theta(z_1,z')d\theta(z_2,z')$. 
Since the map $z\mapsto z'$ is orientation preserving this proves what we want.
\subsection{Step 3} \label{ref-5.3-34} Now we assume $n\ge 3$. First we observe
\label{ref-5.3-35}
\begin{align*}
\beta_{n,m}&=
(-1)^n\frac{1}{(2\pi)^{2n+m-2} }
\int_{B_-\setminus D}
  d\varphi(w,z_2)^md\varphi(z_1,z_2)\cdots d\varphi(z_{n-1},z_n)d\varphi(z_n,z_{1})
d\varphi(w,z_3)\cdots d\varphi(w,z_{n}) \\
&=(-1)^n\frac{1}{(2\pi)^{2n+m-2} }
\int_{\partial(B_-\setminus D)}
 \varphi(w,z_2)^md\varphi(z_1,z_2)\cdots d\varphi(z_{n-1},z_n)d\varphi(z_n,z_{1})
d\varphi(w,z_3)\cdots d\varphi(w,z_{n})
\end{align*}
The components of $\partial(B_-\setminus D)$ which yield non-zero
contributions to the integral have been listed in \S\ref{ref-3.3-8}. Note
that the contribution of \eqref{ref-3.6-10} is zero by \eqref{ref-5.6-33} (integrate
first over $z_i$).
We now list the contributing integrals with the correct signs. 
\begin{enumerate}
\item \label{ref-1-36} Consider the boundary component $\partial_1$ given by contracting
  all edges.  The coordinates in the neighborhood of this component
  are given by $r,x'_2,y'_2,\ldots,x'_n,y'_n$
\begin{align*}
y_1&=1+r\\
x_i&=rx'_i\\
y_i&=1+ry'_i
\end{align*}
with $r\ge 0$. In these coordinate the orientation form on $B_-$ is (up to
a positive multiple) given by
\[
-drdx'_2dy'_2\cdots dx'_ndy'_n
\]
and hence $\partial_1$ is oriented by $dx'_2dy'_2\cdots dx'_ndy'_n$. Whence
the contribution to the integral is equal to $(-1)^n\epsilon_{n,m}$. 
\item Now consider the boundary component $\partial_2$ given by
  contracting $\{z_1,w\}$. We claim $\partial_2\cong C_{n,0}$,
  taking orientations into account.  Coordinates in a neighborhood of
  $\partial_2$ are given by $r,x_2,y_2,\ldots,x_n,y_n$ where $y_1=1+r$
  and $r\ge 0$. Thus the orientation form is $-dr dx_2 dy_2\cdots
  x_n,y_n$ and hence $\partial_2$ is oriented by $dx_2 dy_2\cdots
  x_n,y_n$ which is the normal orientation on $C_{n,0}$ \cite{AMM1}.
  We obtain
\begin{align*}
&(-1)^n\frac{1}{(2\pi)^{2n+m-2} }\int_{\partial_2 }
  \varphi(w,z_2)^md\varphi(z_1,z_2)\cdots d\varphi(z_{n-1},z_n)d\varphi(z_n,z_{1})
d\varphi(w,z_3)\cdots d\varphi(w,z_{n})\\
&= (-1)^n\frac{1}{(2\pi)^{2n+m-2} }\int_{C_{n,0}} \varphi(w,z_2)^md\varphi(w,z_2)\cdots d\varphi(z_{n-1},z_n)d\varphi(z_n,z_1)
d\varphi(w,z_3)\cdots d\varphi(w,z_{n})\\
&=(-1)^n\frac{1}{(m+1)(2\pi)^{2n+m-2} }
\int_{C_{n,0}} d\varphi(w,z_2)^{m+1}d\varphi(z_2,z_3)\cdots d\varphi(z_{n-1},z_n)d\varphi(z_n,z_1)
d\varphi(w,z_3)\cdots d\varphi(w,z_{n})\\
&=-\frac{1}{(m+1)(2\pi)^{2n+m-2} }\int_{C_{n,0}} d\varphi(z_2,z_3)\cdots d\varphi(z_{n-1},z_n)d\varphi(z_n,z_1)d\varphi(w,z_2)^{m+1}
d\varphi(w,z_3)\cdots d\varphi(w,z
_{n})\\
&=-\frac{1}{m+1}\alpha_{n-1,m+1}
\end{align*}
\item 
\label{ref-3-37}
Let $\partial_3$ be obtained by contracting $(z_1,z_2)$.  We
  claim that $\partial_3=-C'_2\times B'_-$ where $B'_-$ is the
  analogue of $B_-$ inside $C_{n,0}$ and where $C'_2$ is the part of
  $C_2$ where the two points are not on the same vertical line. Thus
  $C'_2$ and hence $\partial_3$ has two connected components.

  Coordinates in a neighborhood of $\partial_3$ are given by
  $r,\theta,y_1,x_3,y_3,\ldots, x_n,y_n$ where
\begin{align*}
x_2=&r\cos\theta\\
y_2=&y_1+r\sin\theta
\end{align*}
and $r\ge 0$ and $\theta\not\in \{\pi/2,-\pi/2\}$ (the latter because
of the branch cut involving $\varphi(w,z_2)$). Hence we find that the
orientation form on $B_-$ is up to a positive factor given by 
\[
-dy_1drd\theta dx_3dy_3\cdots dx_ndy_n
\]
Therefore the orientation form on $\partial_3$ (using the outward normal) is
given by 
\[
-dy_1d\theta dx_3dy_2\cdots dx_ndy_n=(-d\theta)\times (-dy_1 dx_3dy_3\cdots dx_ndy_n)
\]
which does indeed represent the orientation on $-C'_2\times B'_-$.

The resulting integral is now
\begin{align*}
&(-1)^n\frac{1}{(2\pi)^{2n+m-2} }\int_{\partial_3}
  \varphi(w,z_2)^md\varphi(z_1,z_2)\cdots d\varphi(z_{n-1},z_n)d\varphi(z_n,z_{1})
d\varphi(w,z_3)\cdots d\varphi(w,z_{n})\\
&=(-1)^{n+1}\frac{1}{2}\frac{1}{(2\pi)^{2n-1} }\int_{B'_-}
  d\varphi(z_2,z_3)\cdots d\varphi(z_{n-1},z_n)d\varphi(z_n,z_{2})
d\varphi(w,z_3)\cdots d\varphi(w,z_{n})\\
&=(-1)^{n+1} \frac{1}{2}\beta_{n-1,1}
\end{align*}
where the factor $1/2$ comes from the fact that $\varphi(w,z_2)=0$ if $z_2$
is to the left of $z_1$. 
\item Let $\partial_4$ be obtained by contracting $(z_2,z_3)$. A
  similar computation as above shows $\partial_4\cong -C_2\times B'_-$. 
The resulting integral is 
\begin{align*}
&(-1)^n\frac{1}{(2\pi)^{2n+m-2} }\int_{\partial_4}
  \varphi(w,z_2)^md\varphi(z_1,z_2)\cdots d\varphi(z_{n-1},z_n)d\varphi(z_n,z_{1})
d\varphi(w,z_3)\cdots d\varphi(w,z_{n})\\
&=(-1)^{n+1}\frac{1}{(2\pi)^{2n+m-2} }\int_{\partial_4}
  \varphi(w,z_2)^md\varphi(z_2,z_3)d\varphi(z_1,z_2)d\varphi(z_3,z_4)\cdots d\varphi(z_{n-1},z_n)\times \\& \qquad \qquad\qquad d\varphi(z_n,z_{1})
d\varphi(w,z_3)\cdots d\varphi(w,z_{n})\\
&=(-1)^{n}\frac{1}{(2\pi)^{2n+m-3} }\int_{B'_-}
  \varphi(w,z_3)^md\varphi(z_1,z_3)d\varphi(z_3,z_4)\cdots d\varphi(z_{n-1},z_n)\times \\& \qquad \qquad\qquad d\varphi(z_n,z_{1})
d\varphi(w,z_3)\cdots d\varphi(w,z_{n})\\
&=(-1)^{n}\frac{1}{(m+1)(2\pi)^{2n+m-3} }\int_{B'_-}
  d\varphi(z_1,z_3)d\varphi(z_3,z_4)\cdots d\varphi(z_{n-1},z_n)\times \\& \qquad \qquad\qquad d\varphi(z_n,z_{1})
d\varphi(w,z_3)^{m+1}\cdots d\varphi(w,z_{n})\\
&=(-1)^{n} \frac{1}{m+1} \beta_{n-1,m+1}
\end{align*}
\end{enumerate}
Combining all contributions we get \eqref{ref-5.3-27}.
\section{A recursion relation for $\bar{\beta}_{n,m}$}
\label{ref-6-38}
 By our definition in \S\ref{ref-3.4-12}, $\bar{\beta}_{n,m}$ is the
integral associated to the enhanced graph 
\[
\psfrag{zn1}[l][l]{$z_{n-1}$}
\psfrag{z1}[][]{$z_1$}
\psfrag{z2}[][]{$z_2$}
\psfrag{zn}[][]{$z_n$}
\psfrag{m}[l][l]{$m$}
\psfrag{w}[l][l]{$w$}
 \raisebox{-2.5cm}{\includegraphics[width=4cm]{wheel_dashed_m_mirror.eps}}
\]
Thus with our standard edge ordering
\[
\bar{\beta}_{n,m}=\frac{1}{(2\pi)^{2n+m-2} }
\int_{B_-\setminus D}
 d\varphi(z_1,z_2)\cdots d\varphi(z_{n-1},z_n)d\varphi(z_n,z_{1})
 d\varphi(w,z_2)\cdots d\varphi(w,z_{n})^m
\]
where
\[
B_-=\{\Re z_1=\Re w+, \Im {z_1}>\Im w\}\subset C_{n+1,0}
\]
\[
D=\{\Re z_n=\Re w, \Im {z_2}>\Im w\}\subset C_{n+1,0}
\]
and where $B_-$ is oriented with the outgoing normal (see \eqref{ref-5.1-25}).
Clearly we have
\[
\bar{\beta}_{2,m}=\beta_{2,m}
\]
\[
\bar{\beta}_{n,1}=\beta_{n,1}
\]
In this section we will prove for $n\ge 3$ 
\begin{equation}
\label{ref-6.1-39}
\bar{\beta}_{n,m}=
(-1)^{n} \frac{1}{2}\beta_{n-1,1}+(-1)^{n+1} \frac{1}{m+1} \bar{\beta}_{n-1,m+1}
+\bar{\epsilon}_{n,m}
\end{equation}
where $\bar{\epsilon}_{n,m}$ is the integral associated to the graph
\[
\psfrag{zn}[][]{$z_n$}
\psfrag{zn1}[l][l]{$z_{n-1}$}
\psfrag{z1}[][]{$z_1$}
\psfrag{z2}[][]{$z_2$}
\psfrag{m}[l][l]{$m$}
\psfrag{w}[l][l]{$w$}
\int \raisebox{-1.5cm}{\includegraphics[width=4cm]{wheel_dashed_bold_m_mirror_boxed.eps}}
\]
(see \S\ref{ref-3.4-12}). Thus with our standard orientation on edges
\[
\bar{\epsilon}_{n,m}=
\frac{1}{(2\pi)^{2n+m-2}}\int \theta(w,z_n)^m d\theta(z_1,z_2)d\theta(z_2,z_3)
\cdots d\theta(z_n,z_1) d\theta(w,z_2)\cdots d\theta(w,z_{n-1})
\]
Here we put $w=0$, $z_1=i$ and the integral is over the complement
of the diagonal in $(\CC-\{0,i\})^{n-1}$.
\subsection{Proof}
\label{ref-6.1-40}
We assume  $n\ge 3$. We first observe 
\begin{align*}
\bar{\beta}_{n,m}&=
\frac{1}{(2\pi)^{2n+m-2} }
\int_{B_-\setminus D}
  d\varphi(w,z_n)^md\varphi(z_1,z_2)\cdots d\varphi(z_{n-1},z_n)d\varphi(z_n,z_{1})
d\varphi(w,z_2)\cdots d\varphi(w,z_{n-1}) \\
&=\frac{1}{(2\pi)^{2n+m-2} }
\int_{\partial(B_-\setminus D)}
 \varphi(w,z_n)^md\varphi(z_1,z_2)\cdots d\varphi(z_{n-1},z_n)d\varphi(z_n,z_{1})
d\varphi(w,z_2)\cdots d\varphi(w,z_{n-1})
\end{align*}
The components of $\partial(B_-\setminus D)$ which yield non-zero
contributions to the integral have been listed in \S\ref{ref-3.3-8}. Note
that the contribution of \eqref{ref-3.7-11} is zero by \eqref{ref-5.6-33} (integrate
first over $z_{n-1}$).
We now list the contributing integrals with the correct signs. 
\begin{enumerate}
\item Consider the boundary component $\partial_1$ given by contracting
  all edges. Then a computation similar to \S\ref{ref-5.3-34}(\ref{ref-1-36})
yields that the contribution is $\bar{\epsilon}_{n,m}$. 
\item Now consider the boundary component $\partial_2$ given by 
contracting $(z_{n-1},w)$.  It this case we get an integral which contains
as a subintegral (up to a scalar factor)
\[
\int_{z_n\in \Hscr\setminus\{w,z_1\}} d\varphi(w,z_n)^{m+1} d\varphi(z_n,z_1) 
\]
Using \eqref{ref-4.3-18} this is equal to
\[
(2\pi)^{m+1}\varphi(w,z_1)-2\pi\varphi(w,z_1)^{m+1}=(2\pi)^{m+1} \times 2\pi-
2\pi\times (2\pi)^{m+1}=0
\]
So we are lucky that there is no contribution in this case.
\item Let $\partial_3$ be obtained by contracting $(z_1,z_n)$. A
  similar computation as in \S\ref{ref-5.3-34}\ref{ref-3-37}) shows
  $\partial_3=-C'_2\times B'_-$. The resulting integral is
\begin{align*}
&\frac{1}{(2\pi)^{2n+m-2} }\int_{\partial_3}
  \varphi(w,z_n)^md\varphi(z_1,z_2)\cdots d\varphi(z_{n-1},z_n)d\varphi(z_n,z_{1})
d\varphi(w,z_2)\cdots d\varphi(w,z_{n-1})\\
&=(-1)^{n-1}\frac{1}{(2\pi)^{2n+m-2} }\int_{\partial_3}
 \varphi(w,z_n)^m d\varphi(z_n,z_{1})d\varphi(z_1,z_2)\cdots d\varphi(z_{n-1},z_n)
d\varphi(w,z_2)\cdots d\varphi(w,z_{n-1})\\
&=(-1)^{n}\frac{1}{2}\frac{1}{(2\pi)^{2n-2} }\int_{B'_-}
d\varphi(z_1,z_2)\cdots d\varphi(z_{n-1},z_1)
d\varphi(w,z_2)\cdots d\varphi(w,z_{n-1})\\
&=(-1)^{n} \frac{1}{2}\bar{\beta}_{n-1,1}=(-1)^{n}\frac{1}{2}\beta_{n-1,1}
\end{align*}
where the factor $1/2$ comes from the fact that $\varphi(w,z_n)=0$ if $z_2$
is to the left of $z_1$. 
\item Let $\partial_4$ be obtained by contracting $(z_{n-1},z_n)$.  A
  similar computation as in \S\ref{ref-5.3-34}\ref{ref-3-37}) shows
  $\partial_4\cong -C_2\times B'_-$. The resulting integral is
\begin{align*}
&\frac{1}{(2\pi)^{2n+m-2} }\int_{\partial_4}
  \varphi(w,z_n)^md\varphi(z_1,z_2)\cdots d\varphi(z_{n-1},z_n)d\varphi(z_n,z_{1})
d\varphi(w,z_2)\cdots d\varphi(w,z_{n-1})\\
&=(-1)^{n}\frac{1}{(2\pi)^{2n+m-2} }\int_{\partial_4}
  \varphi(w,z_n)^md\varphi(z_{n-1},z_n)d\varphi(z_1,z_2)\cdots d\varphi(z_{n-2},z_{n-1})\times \\& \qquad \qquad\qquad d\varphi(z_n,z_{1})
d\varphi(w,z_2)\cdots d\varphi(w,z_{n-1})\\
&=(-1)^{n+1}\frac{1}{(2\pi)^{2n+m-3} }\int_{B'_-}
  \varphi(w,z_{n-1})^md\varphi(z_1,z_2)\cdots d\varphi(z_{n-2},z_{n-1})\times \\& \qquad \qquad\qquad d\varphi(z_{n-1},z_{1})
d\varphi(w,z_2)\cdots d\varphi(w,z_{n-1})\\
&=(-1)^{n+1}\frac{1}{(m+1)(2\pi)^{2n+m-3} }\int_{B'_-}
  d\varphi(z_1,z_2)\cdots d\varphi(z_{n-2},z_{n-1})\times \\& \qquad \qquad\qquad d\varphi(z_{n-1},z_{1})
d\varphi(w,z_2)\cdots d\varphi(w,z_{n-1})^{m+1}\\
&=(-1)^{n+1} \frac{1}{m+1} \bar{\beta}_{n-1,m+1}
\end{align*}
Combining all contributions we get \eqref{ref-6.1-39}.
\end{enumerate}
\section{A recursion not involving unknown quantities.}
\label{ref-7-41}
Put
\[
\hat{\beta}_{n,m}=\frac{1}{2}(\beta_{n,m}+(-1)^{n}\bar{\beta}_{n,m})
\]
We claim
\begin{equation}
\label{ref-7.1-42}
\hat{\beta}_{n,m}=-\frac{1}{2(m+1)}\alpha_{n-1,m+1}
+(-1)^{n+1} \frac{1}{2}\hat{\beta}_{n-1,1}+(-1)^n\frac{1}{m+1} \hat{\beta}_{n-1,m+1}
\end{equation}
with initial condition
\begin{equation}
\label{ref-7.2-43}
\hat{\beta}_{2,m}=-\frac{1}{8}+\frac{1}{m+1}\left( \frac{1}{2}-\frac{1}{m+2}\right)
\end{equation}
\subsection{Step 1} 
\label{ref-7.1-44}
We claim first that 
\begin{equation}
\label{ref-7.3-45}
\epsilon_{n,m}=-\bar{\epsilon}_{n,m}
\end{equation}
We have
\[
\epsilon_{n,m}=
\frac{1}{(2\pi)^{2n+m-2}}\int \theta(w,z_2)^m d\theta(z_1,z_2)d\theta(z_2,z_3)
\cdots d\theta(z_n,z_1) d\theta(w,z_3)\cdots d\theta(w,z_n)
\]
Here we put $w=0$, $z_1=i$ and the integral is over a suitable open
subset of $\CC^{2n-1}$. 

We apply the   permutation $(2n)(3\,
n-1)\cdots $ to $z_1,\ldots,z_n$. This does not change the orientation. We find
\begin{align*}
\epsilon_{n,m}&=
\frac{1}{(2\pi)^{2n+m-2}}\int \theta(w,z_n)^m d\theta(z_1,z_{n})d\theta(z_n,z_{n-1})
\cdots d\theta(z_2,z_1) d\theta(w,z_{n-1})\cdots d\theta(w,z_2)\\
&=\frac{1}{(2\pi)^{2n+m-2}}\int \theta(w,z_n)^m d\theta(z_n,z_{1})d\theta(z_{n-1},z_{n})
\cdots d\theta(z_1,z_2) d\theta(w,z_{n-1})\cdots d\theta(w,z_2)\\
&=(-1)^s \bar{\epsilon}_{n,m}
\end{align*}
The sign $(-1)^s$ comes from the fact that we have invert the order of the
1-forms in $d\theta(z_n,z_1)d\theta(z_{n-1} z_{n})\cdots d\theta(z_1,z_2)$ and
in $d\theta(w,z_{n-1})\cdots d\theta(w,z_2)$.
Thus
\[
s=\frac{n(n-1)}{2}+\frac{(n-2)(n-3)}{2} =n^2-3n+3
\]
which is always an odd number. This proves \eqref{ref-7.3-45}.
\subsection{Step 2}
\label{ref-7.2-46}
Formula \eqref{ref-7.2-43} is obvious from \eqref{ref-5.2-26}
and the fact that $\hat{\beta}_{2,m}=\beta_{2,m}$. Formula \eqref{ref-7.1-42} follows from \eqref{ref-5.3-27} and \eqref{ref-6.1-39} together
with the following computation
\begin{align*}
\hat{\beta}_{n,m}&=\frac{1}{2}\left(-\frac{1}{m+1}\alpha_{n-1,m+1}
+(-1)^{n+1} \frac{1}{2}\beta_{n-1,1}+\frac{1}{2}\beta_{n-1,1}
+(-1)^{n} \frac{1}{m+1} \beta_{n-1,m+1}-\frac{1}{m+1} \bar{\beta}_{n-1,m+1}\right)\\
&=-\frac{1}{2(m+1)}\alpha_{n-1,m+1}
+(-1)^{n+1} \frac{1}{2}\hat{\beta}_{n-1,1}+(-1)^n\frac{1}{m+1} \hat{\beta}_{n-1,m+1}
\end{align*}
\section{A formula for $w_n$}
\label{ref-8-47}
In this section we prove the following fact
\begin{equation}
\label{ref-8.1-48}
w_n=
\begin{cases}
\hat{\beta}_{n,1}-\frac{1}{2}\alpha_{n-1,2}&\text{if $n$ is even}\\
0&\text{if $n$ is odd}
\end{cases}
\end{equation}
\subsection{Step 1} We assume $n$ odd. Put $w=i$. We identify
\label{ref-8.1-49}
$C_{n,0}$ with an open subset of $\Hscr^{n}$. The orientation is
derived from the standard orientation on $\Hscr^n$ (according to
\cite{AMM1}). 

We consider the map $\alpha:\Hscr\r \Hscr: x+iy\mapsto -x+iy$ and we extend
$\alpha$ diagonally to a map $\Hscr^n\r \Hscr^n$ also denoted by $\alpha$.
$\alpha$ multiplies the orientation by $(-1)^n$.

Hence we have
\[
w_n=(-1)^n\frac{1}{(2\pi)^{2n}} \int_{C_{n,0}}
d\varphi(\alpha(z_1),\alpha(z_2))\cdots d\varphi(\alpha(z_1),\alpha(z_n))
d\varphi(w,\alpha(z_1))\cdots d\varphi(w,\alpha(z_n))
\]
(since $\alpha(w)=w$).

It is now easy to see that $\varphi(\alpha(u),\alpha(v))=2\pi-\varphi(u,v)$ and
hence $d\varphi(\alpha(u),\alpha(v))=-d\varphi(u,v)$. Thus
\[
w_n=(-1)^n(-1)^{2n}w_n=(-1)^nw_n
\]
Since $n$ is odd this implies $w_n=0$.  
\subsection{Step 2} Consider the following enhanced graph
\label{ref-8.2-50}
\begin{equation}
\label{ref-8.2-51}
\psfrag{z1}[][]{$z_1$}
\psfrag{z2}[][]{$z_2$}
\psfrag{zn}[][]{$z_n$}
\psfrag{m}[][]{$$}
\psfrag{w}[][]{$w$}
\raisebox{-2cm}{\includegraphics[width=4cm]{wheel_n_bold_euclidean.eps}}
\end{equation}
We claim that its corresponding integral  
\[
\delta=\frac{1}{(2\pi)^{2n}} \int_{C_{n+1}} \theta(w,z_1) d\theta(z_1,z_2)\cdots
d\theta(z_n,z_1) d\theta(w,z_2)\cdots d\theta(w,z_n)
\]
is zero when $n$ is even.

It's easy to see that permuting the points does not change the
orientation on $C_{n+1}$.  We apply the permutation $(2n)(3\,
n-1)\cdots $ for $z_1,z_2,\ldots,z_n$. Hence we also have
\begin{align*}
\delta&=\frac{1}{(2\pi)^{2n}} \int_{C_{n+1}}
\theta(w,z_1) d\theta(z_1,z_n)d\theta(z_n z_{n-1})\cdots d\theta(z_2,z_1)
d\theta(w,z_{n})\cdots d\theta(w,z_2)\\
&=\frac{1}{(2\pi)^{2n}} \int_{C_{n+1}}
\theta(w,z_1) d\theta(z_n,z_1)d\theta(z_{n-1} z_{n})\cdots d\theta(z_1,z_2)
d\theta(w,z_{n})\cdots d\theta(w,z_2)
\end{align*}
To put this back in standard form we have to invert the order of the
1-forms in $d\theta(z_n,z_1)d\theta(z_{n-1} z_{n})\cdots d\theta(z_1,z_2)$ and
in $d\theta(w,z_{n})\cdots d\theta(w,z_2)$. This introduces a cumulative sign
of $(-1)^s$ where 
\[
s=\frac{n(n-1)}{2}+\frac{(n-1)(n-2)}{2} =(n-1)^2
\]
Thus 
\[
\delta=-(-1)^n \delta
\]
and hence if $n$ is even then $\delta=0$. 
\subsection{Step 3} Now we assume $n$ even. 
\label{ref-8.3-52}
As a pedagogical device we keep signs of the form $(-1)^n$ below. 
 We first observe
\begin{align*}
w_n&=(-1)^n \frac{1}{(2\pi)^{2n}}
\int_{C_{n+1,0}}
 d\varphi(w,z_1)  d\varphi(z_1,z_2)\cdots d\varphi(z_{n-1},z_n)d\varphi(z_n,z_{1}) d\varphi(w,z_2)\cdots
d\varphi(w,z_n)\\
&=(-1)^n \frac{1}{(2\pi)^{2n}}\int_{\text{boundary}}
 \varphi(w,z_1)  d\varphi(z_1,z_2)\cdots d\varphi(z_{n-1},z_n)d\varphi(z_n,z_{1}) d\varphi(w,z_2)\cdots
d\varphi(w,z_n)
\end{align*}
The boundary consists of two parts. One part consists of both sides
of the branch cut $B_{\pm}$.
\[
B_-=\{\Re z_1=\Re w+, \Im {z_1}>\Im w\}\subset C_{n+1,0}
\]
\[
B_+=\{\Re z_1=\Re w-, \Im {z_1}>\Im w\}\subset C_{n+1,0}
\]
The values of $\varphi(w,-)$ on both sides of the branch cut differ by
$2\pi$. Hence the contribution of this  part of the boundary is equal to
\[
(-1)^n \frac{1}{(2\pi)^{2n-1}}\int_{B_-}
 \varphi(w,z_1)  d\varphi(z_1,z_2)\cdots d\varphi(z_{n-1},z_n)d\varphi(z_n,z_{1}) d\varphi(w,z_2)\cdots
d\varphi(w,z_n)
\]
which is nothing but
\[
(-1)^n \beta_{n,1}
\]
In \S\ref{ref-3.2-5} we have listed the ``classical'' parts of the boundary which
contribute to the integral. Here we give the contributions with the precise 
signs. Note that the contribution of \eqref{ref-8.2-51} is zero
by Step 2.
\begin{enumerate}
\item Let $\partial_1$ be the boundary component given by contracting
  $(w,z_n)$. According to \cite{AMM1} we have $\partial_1\cong-C_2\times C_{n,0}$. 
The contribution to the integral is
\begin{multline*}
(-1)^n \frac{1}{(2\pi)^{2n}}\int_{\partial_1}
 \varphi(w,z_1)  d\varphi(z_1,z_2)\cdots d\varphi(z_{n-1},z_n)d\varphi(z_n,z_{1}) d\varphi(w,z_2)\cdots
d\varphi(w,z_{n-1})d\varphi(w,z_n)\\
=(-1)^n \frac{1}{(2\pi)^{2n}}\int_{\partial_1}
 \varphi(w,z_1)  d\varphi(w,z_n) d\varphi(z_1,z_2)\cdots d\varphi(z_{n-1},z_n)d\varphi(z_n,z_{1}) d\varphi(w,z_2)\cdots
d\varphi(w,z_{n-1})\\
=(-1)^{n+1} \frac{1}{(2\pi)^{2n-1}}\int_{C_{n,0}}
 \varphi(w,z_1)  d\varphi(z_1,z_2)\cdots d\varphi(z_{n-1},w)d\varphi(w,z_{1}) d\varphi(w,z_2)\cdots
d\varphi(w,z_{n-1})\\
=(-1)^{n+1} \frac{1}{2(2\pi)^{2n-1}}
\int_{C_{n,0}}
 d\varphi(z_1,z_2)\cdots d\varphi(z_{n-1},w)d\varphi(w,z_{1})^{2} d\varphi(w,z_2)\cdots
d\varphi(w,z_{n-1})\\
=(-1)^{n+1} \frac{1}{2}\alpha_{n-1,2}
\end{multline*}
\item Let $\partial_2\cong -C_2\times C_{n,0}$ be obtained by
  contracting $(z_1,z_2)$. We obtain
\begin{multline*}
(-1)^n \frac{1}{(2\pi)^{2n}}\int_{\partial_2} 
 \varphi(w,z_1)  d\varphi(z_1,z_2)\cdots d\varphi(z_{n-1},z_n)d\varphi(z_n,z_{1}) d\varphi(w,z_2)\cdots
d\varphi(w,z_n)\\
=(-1)^{n+1} \frac{1}{(2\pi)^{2n-1}}\int_{C_{n,0}}
 \varphi(w,z_2)  d\varphi(z_2,z_3)\cdots d\varphi(z_{n-1},z_n)d\varphi(z_n,z_{2}) d\varphi(w,z_2)\cdots
d\varphi(w,z_n)\\
=(-1)^{n+1} \frac{1}{(2)(2\pi)^{2n-1}}
\int_{C_{n,0}} 
d\varphi(z_2,z_3)\cdots d\varphi(z_{n-1},z_n)d\varphi(z_n,z_{2}) d\varphi(w,z_2)^{2}\cdots
d\varphi(w,z_n)
%\\
%(-1)^{n+1}\frac{1}{2}\gamma_{n-1,2}
\end{multline*}
\item Let $\partial_3\cong -C_2\times C_{n,0}$ be obtained by
  contracting $(z_n,z_1)$. We obtain
\begin{multline*}
(-1)^n \frac{1}{(2\pi)^{2n}}\int_{\partial_3} 
 \varphi(w,z_1)  d\varphi(z_1,z_2)\cdots d\varphi(z_{n-1},z_n)d\varphi(z_n,z_{1}) d\varphi(w,z_2)\cdots
d\varphi(w,z_n)\\
=- \frac{1}{(2\pi)^{2n}}\int_{\partial_3} 
 \varphi(w,z_1) d\varphi(z_n,z_{1})  d\varphi(z_1,z_2)\cdots d\varphi(z_{n-1},z_{n})  d\varphi(w,z_2)\cdots
d\varphi(w,z_n)\\
=\frac{1}{(2\pi)^{2n-1}}\int_{C_{n,0}} 
 \varphi(w,z_1)   d\varphi(z_1,z_2)\cdots d\varphi(z_{n-1},z_{1})  d\varphi(w,z_2)\cdots
d\varphi(w,z_{n-1}) d\varphi(w,z_1)\\
=(-1)^{n-2}\frac{1}{(2\pi)^{2n-1}}\int_{C_{n,0}} 
 \varphi(w,z_1)   d\varphi(z_1,z_2)\cdots d\varphi(z_{n-1},z_{1})  d\varphi(w,z_1) d\varphi(w,z_2)\cdots
d\varphi(w,z_{n-1})\\
=(-1)^{n-2}\frac{1}{(2)(2\pi)^{2n-1}}\int_{C_{n,0}} 
   d\varphi(z_1,z_2)\cdots d\varphi(z_{n-1},z_{1})  d\varphi(w,z_1)^{2} d\varphi(w,z_2)\cdots
d\varphi(w,z_{n-1})%\\
%=(-1)^{n-2} \frac{1}{m+1}\gamma_{n-1,m+1}
\end{multline*}
Hence the contributions of $\partial_2$ and $\partial_3$ cancel. The
contributions of $\partial_1$ and $B_-$ yield \eqref{ref-8.1-48} in case $n$
is even.
\end{enumerate}
\section{Solving the recursion}
\label{ref-9-53}
\subsection{The problem}
\label{ref-9.1-54}
For the benefit of the reader we restate the recursion relations we have
derived. See \eqref{ref-4.1-15}\eqref{ref-4.2-16}\eqref{ref-7.1-42}\eqref{ref-7.2-43}\eqref{ref-8.1-48}.
\begin{align*}
\alpha_{0,m}&= 1\\
\alpha_{n,m}&=(-1)^n \left( \frac{1}{2}\alpha_{n-1,2}-\frac{1}{m+1} \alpha_{n-1,m+1}\right) \qquad \text{(for $n\ge 1$)}\\
\hat{\beta}_{2,m}&=-\frac{1}{8}+\frac{1}{m+1}\left( \frac{1}{2}-\frac{1}{m+2}\right)\\
\hat{\beta}_{n,m}&=-\frac{1}{2(m+1)}\alpha_{n-1,m+1} +(-1)^{n+1}
\frac{1}{2}\hat{\beta}_{n-1,1}+(-1)^n\frac{1}{m+1}
\hat{\beta}_{n-1,m+1 }
\qquad \text{(for $n\ge 3$)}
\\
w_n&= \hat{\beta}_{n,1}-\frac{1}{2}\alpha_{n-1,2}\qquad\text{(for $n\ge 2$ even)}
\end{align*}
\subsection{Eliminating some signs and fractions}
\label{ref-9.2-55}
We simplify  the equations by putting
\[
\tilde{\alpha}_{n,m}=(-1)^{\frac{n(n+1)}{2}}\frac{1}{m!} \alpha_{n,m}
\]
Then we have
\[
\tilde{\alpha}_{0,m}=\frac{1}{m!}
\]
and
\begin{equation}
\label{ref-9.1-56}
\tilde{\alpha}_{n,m}=  
\frac{1}{m!}\tilde{\alpha}_{n-1,2}- \tilde{\alpha}_{n-1,m+1}
\end{equation}
Similarly we put
\[
\tilde{\beta}_{n,m}=(-1)^{\frac{n(n+1)}{2}} \frac{\hat{\beta}_{n,m}}{m!}
\]
so that we get. 
\[
\tilde{\beta}_{2,m}=\frac{1}{8m!}-\frac{1}{2(m+1)!}+\frac{1}{(m+2)!}
\]
and
\begin{equation}
\label{ref-9.2-57}
\tilde{\beta}_{n,m}=-(-1)^{n}\frac{1}{2}\tilde{\alpha}_{n-1,m+1}-\frac{1}{2m!}\tilde{\beta}_{n-1,1}+\tilde{\beta}_{n-1,m+1}
\end{equation}
\subsection{Computing $\tilde{\alpha}$}
\label{ref-9.3-58}
 Iterating \eqref{ref-9.1-56} we find
\begin{align*}
\tilde{\alpha}_{n,2}&=  \frac{1}{2!}\tilde{\alpha}_{n-1,2}
- \tilde{\alpha}_{n-1,3}\\
&=\frac{1}{2!}\tilde{\alpha}_{n-1,2}
-\frac{1}{3!} \tilde{\alpha}_{n-2,2}+\tilde{\alpha}_{n-2,4}\\
&=\left(\frac{1}{2!}\tilde{\alpha}_{n-1,2}
-\frac{1}{3!} \tilde{\alpha}_{n-2,2}+\cdots + (-1)^{n+1}\frac{1}{(n+1)!} \tilde{\alpha}_{0,2}\right)+(-1)^n\tilde{\alpha}_{0,n+2}
\end{align*}
Thus if we put 
\[
A_2=\sum_{n\ge 0} \tilde{\alpha}_{n,2} x^n
\]
we get
\[
A_2=\left(\frac{x}{2!}-\frac{x^2}{3!}+\cdots\right)A_2+\sum_{n\ge 0}(-1)^n \frac{x^n}{(n+2)!}
\]
Here
\[
\frac{x}{2!}-\frac{x^2}{3!}+\cdots=\frac{e^{-x}-1+x}{x}
\]
and
\[
\sum_{n\ge 0}(-1)^n \frac{x^n}{(n+2)!}=\frac{e^{-x}-1+x}{x^2}
\]
from which we deduce
\[
A_2=-\frac{1}{x}+\frac{1}{1-e^{-x}}
\]
For use below we record
\[
A'_2\overset{\text{def}}{=}A_2-A_2(0)=-\frac{1}{x}+\frac{1+e^{-x}}{2(1-e^{-x})}
\]

Knowing $\tilde{\alpha}_{n,2}$ we can compute arbitrary $\tilde{\alpha}_{n,m}$ via
\begin{equation}
\label{ref-9.3-59}
\begin{split}
\tilde{\alpha}_{n,m}&=\frac{1}{m!}\tilde{\alpha}_{n-1,2}- \tilde{\alpha}_{n-1,m+1}\\
&=\frac{1}{m!}\tilde{\alpha}_{n-1,2}-\frac{1}{(m+1)!}\tilde{\alpha}_{n-2,2}
+\cdots+(-1)^{n-1}\frac{1}{(n+m-1)!}\tilde{\alpha}_{0,2}
+
(-1)^n\frac{1}{(m+n)!}
\end{split}
\end{equation}
We will not write down a closed expression for $\tilde{\alpha}_{n,m}$
as we won't need it.

\subsection{Computing a sum}
\label{ref-9.4-60}
For use below we compute the sum for $n\ge 2$
\[
s_n=(-1)^{n-1}\tilde{\alpha}_{n-1,2}+(-1)^{n-2}\tilde{\alpha}_{n-2,3}
+\cdots +\tilde{\alpha}_{2,n-1}
\]
where an empty sum is interpreted as $0$. Using \eqref{ref-9.3-59} we get
\[
\tilde{\alpha}_{2,n-1}=\frac{1}{(n-1)!}\tilde{\alpha}_{1,2}-\frac{1}{n!}\tilde{\alpha}_{0,2}+\frac{1}{(n+1)!}
\]
\[
\tilde{\alpha}_{3,n-2}=
\frac{1}{(n-2)!} \tilde{\alpha}_{2,2}-
\frac{1}{(n-1)!} \tilde{\alpha}_{1,2}+
\frac{1}{n!} \tilde{\alpha}_{0,2}-
\frac{1}{(n+1)!} 
\]
and
\[
\tilde{\alpha}_{n-1,2}=\frac{1}{2!} \tilde{\alpha}_{n-2,2}-\frac{1}{3!}
\tilde{\alpha}_{n-3,2}+\cdots + (-1)^{n-2} \frac{1}{n!}\tilde{\alpha}_{0,2}
+(-1)^{n-1} \frac{1}{(n+1)!}
\]
so that we get
\[
s_n=
(-1)^{n-1}\frac{1}{2!} \tilde{\alpha}_{n-2,2}+
(-1)^{n-2} \frac{2}{3!} \tilde{\alpha}_{n-3,2}+\cdots+
(n-2)\frac{1}{(n-1)!}\tilde{\alpha}_{1,2}-(n-2)\frac{1}{n!}\tilde{\alpha}_{0,2}+ (n-2)\frac{1}{(n+1)!}
\]
Put
\[
T=\frac{x^2}{2!}-\frac{2x^3}{3!}+\cdots=-xe^{-x}-e^{-x}+1
\]
Then we get 
\[
s_n=(-1)^{n-1}(TA_2)[x^n]+\frac{1}{n!}\tilde{\alpha}_{0,2}+\frac{n-2}{(n+1)!}
\]
where $(-)[x^n]$ denotes the coefficient of $x^n$.
 We will regard this formula  as a definition if $n=0,1$.  We
obtain
\[
\sum_{n\ge 0} (-1)^n s_n x^n=-TA_2+ \sum_{n\ge 0} (-1)^n \tilde{\alpha}_{0,2}\frac{x^n}{n!}+\sum_{n\ge 0} (-1)^n\frac{n-2}{(n+1)!}x^n
\]
We have
\[
\sum_{n\ge 0} (-1)^n \tilde{\alpha}_{0,2}\frac{x^n}{n!}=\frac{1}{2}e^{-x}
\]
and
\begin{align*}
\sum_{n\ge 0} (-1)^n\frac{n-2}{(n+1)!}x^n&=\sum_{n\ge 0}(-1)^n \frac{n+1}{(n+1)!}x^n-3
\sum_{n\ge 0}(-1)^n \frac{1}{(n+1)!}x^n\\
&=e^{-x}+3\frac{e^{-x}-1}{x}
\end{align*}
A computation with a computer algebra package now yields
\[
\sum_{n\ge 0} (-1)^ns_n x^n=
\frac{-(x+4)e^{-2x}+(2x^2+3x+8)e^{-x}-(2x+4)}{2(1-e^{-x})x}
\]
Hence
\[
S\overset{\text{def}}{=}\sum_{n\ge 0} s_n x^n=
\frac{(x-4)e^{2x}+(2x^2-3x+8)e^{x}+(2x-4)}{2(e^{x}-1)x}
\]
We have $S=-3/2+O(x^4)$. For use below we record
\[
S'\overset{\text{def}}{=}\sum_{n\ge 2} s_n x^n=S+\frac{3}{2}=
\frac{(x-4)e^{2x}+(2x^2+8)e^x-x-4}{2(e^x-1)x}
\]
\subsection{Computing $\tilde{\beta}_{n,1}$}
\label{ref-9.5-61}
We will only compute $\tilde{\beta}_{n,1}$ since this the only thing we need.
Iterating \eqref{ref-9.2-57} we get for $n\ge 3$
\begin{align*}
  \tilde{\beta}_{n,1}&=(-1)^{n-1}\frac{1}{2}\tilde{\alpha}_{n-1,2}-\frac{1}{2\cdot 1!}
\tilde{\beta}_{n-1,1}+\tilde{\beta}_{n-1,2}\\
&=\frac{1}{2}\biggl((-1)^{n-1}\tilde{\alpha}_{n-1,2}+(-1)^{n-2}\tilde{\alpha}_{n-2,3}
+\cdots +\tilde{\alpha}_{2,n-1}\biggr)-\\
&\qquad\qquad\frac{1}{2}\left(\frac{1}{1!}\tilde{\beta}_{n-1,1}
+\frac{1}{2!}\tilde{\beta}_{n-2,1}+\cdots+\frac{1}{(n-2)!}\tilde{\beta}_{2,1}\right)+\tilde{\beta}_{2,n-1} \\
&=\frac{s_n}{2} -\frac{1}{2}\left(\frac{1}{1!}\tilde{\beta}_{n-1,1}
+\frac{1}{2!}\tilde{\beta}_{n-2,1}+\cdots+\frac{1}{(n-2)!}\tilde{\beta}_{2,1}\right)+\tilde{\beta}_{2,n-1}
\end{align*}
Since $s_2=0$ this identity holds for $n=2$ if we interpret an empty
sum as zero. 

If we put
\[
B_1=\sum_{n\ge 2} \tilde{\beta}_{n,1}x^n
\]
then we have
\[
B_1=\frac{1}{2} S' -\frac{1}{2}\left( \frac{x}{1!}+\frac{x^2}{2!}+
\cdots\right)B_1
+
\sum_{n\ge 2}x^n\left(\frac{1}{8(n-1)!}-\frac{1}{2n!}+\frac{1}{(n+1)!}\right)
\]
We have
\[
 \frac{x}{1!}+\frac{x^2}{2!}+
\cdots=e^x-1
\]
\[
\sum_{n\ge 2} \frac{x^n}{(n-1)!}=x(e^x-1)
\]
\[
\sum_{n\ge 2} \frac{x^n}{n!}=e^x-1-x
\]
\[
\sum_{n\ge 2} \frac{x^n}{(n+1)!}=\frac{1}{x}\left(e^x-1-x-\frac{x^2}{2} \right)
\]
and hence
\begin{align*}
\sum_{n\ge
  2}x^n\left(\frac{1}{8(n-1)!}-\frac{1}{2n!}+\frac{1}{(n+1)!}\right)& =
\frac{1}{8} x(e^x-1)-\frac{1}{2}
(e^x-1-x)+\frac{1}{x}\left(e^x-1-x-\frac{x^2}{2}\right)\\
&=\left(\frac{x}{8}-\frac{1}{2}+\frac{1}{x}\right)e^x-\frac{1}{x}
-\frac{1}{2}-\frac{x}{8}
\end{align*}
Thus we get
\[
\left(1+\frac{1}{2}(e^x-1)\right)B_1=\frac{1}{2}S'+\left(\frac{x}{8}-\frac{1}{2}+\frac{1}{x}\right)e^x-\frac{1}{x}
-\frac{1}{2}-\frac{x}{8}
\]
Invoking once again a computer algebra package we find 
\[
B_1=\frac{(x-2)e^x+x+2}{4(e^x-1)}
\]
\subsection{Computing $w_n$ }
\label{ref-9.6-62}
If $n$ is even we have
\[
(-1)^{\frac{n(n-1)}{2}}w_n=\tilde{\beta}_{n,1}-\tilde{\alpha}_{n-1,2}
\]
Since $A'_2$, $B_1$ are respectively an odd and an even function of
$x$, $B_1(0)=0$  and $w_n=0$ for $n$ odd we obtain
\[
\sum_{n\ge 2}(-1)^{\frac{n(n-1)}{2}}w_nx^n=B_1-xA'_2
\]
and using a computer algebra package we find
\begin{equation}
\label{ref-9.4-63}
\sum_{n\ge 2}(-1)^{\frac{n(n-1)}{2}}w_nx^n=\frac{(x+2)e^{-x}+x-2}{4(e^{-x}-1)}
\end{equation}
The derivative of 
\[
\frac{1}{2}\log \frac{e^{x/2}-e^{-x/2}}{x}
\]
is equal to
\[
-\frac{(x+2)e^{-x}+x-2}{4(e^{-x}-1)x}
\]
Dividing \eqref{ref-9.4-63} by $x$ and integrating finishes the proof
of \eqref{ref-1.1-1}.
%\bibliography{mybibs.bib}

\begin{thebibliography}{10}

\bibitem{AMM1}
D.~Arnal, D.~Manchon, and M.~Masmoudi, {\em Choix des signes pour la
  formalit\'e de {M}. {K}ontsevich}, Pacific J. Math. {\bf 203} (2002), no.~1,
  23--66.

\bibitem{vdbcalaque}
D.~Calaque and M.~Van~den Bergh, {\em Hochschild cohomology and {A}tiyah
  classes}, submitted, arXiv:0708.2725.

\bibitem{CKTB}
A.~Cattaneo, B.~Keller, C.~Torossian, and A.~Brugui{\`e}res, {\em
  Introduction}, D\'eformation, quantification, th\'eorie de Lie (Paris),
  Panor. Synth\`eses, vol.~20, Soc. Math. France, Paris, 2005, Dual
  French-English text, pp.~1--9, 11--18.

\bibitem{CF2}
A.~S. Cattaneo and G.~Felder, {\em Coisotropic submanifolds in {P}oisson
  geometry and branes in the {P}oisson sigma model}, Lett. Math. Phys. {\bf 69}
  (2004), 157--175.

\bibitem{CF3}
\bysame, {\em Relative formality theorem and quantisation of coisotropic
  submanifolds}, Adv. Math. {\bf 208} (2007), no.~2, 521--548.

\bibitem{GSL}
{GSL - GNU Scientific Library}, http://www.gnu.org/software/gsl.

\bibitem{Ko3}
M.~Kontsevich, {\em Deformation quantization of {P}oisson manifolds}, Lett.
  Math. Phys. {\bf 66} (2003), no.~3, 157--216.

\bibitem{Maxima}
Maxima, http://maxima.sourceforge.net.

\bibitem{SAGE}
{SAGE Mathematical Software, Version 2.6}, http://www.sagemath.org.

\bibitem{Shoikhet1}
B.~Shoikhet, {\em On the {D}uflo formula for {$L_\infty$}-algebras and
  {Q}-manifolds}, math/9812009.

\bibitem{Shoikhet}
\bysame, {\em Vanishing of the {K}ontsevich integrals of the wheels}, Lett.
  Math. Phys. {\bf 56} (2001), no.~2, 141--149, EuroConf\'erence Mosh\'e Flato
  2000, Part II (Dijon).

\bibitem{Wolfram}
E.~W. Weisstein, {\em Modified {B}ernoulli {N}umber}, MathWorld,
  http://mathworld.wolfram.com /ModifiedBernoulliNumber.html.


\bibitem{Willwacher} T.~Willwacher  {\em A counterexample to the quantizability of modules}, Lett. Math. Phys. {\bf 81} (2007), no.~3, pp. 265--280.


\end{thebibliography}
%\bibliographystyle{amsabbrv}

\def\cprime{$'$} \def\cprime{$'$} \def\cprime{$'$}
\ifx\undefined\bysame
\newcommand{\bysame}{\leavevmode\hbox to3em{\hrulefill}\,}
\fi

\end{document}